\documentclass[a4paper]{article}

\usepackage[bookmarks=true,linkcolor=black,citecolor=black,urlcolor=black,colorlinks=true,breaklinks=true,bookmarksopen=true]{hyperref}
\usepackage{amsmath}
\usepackage{amsfonts}
\usepackage{nicefrac}
\usepackage{graphicx}
\usepackage{epstopdf}
\usepackage{algorithmic}
\usepackage{enumerate}
\usepackage{wrapfig}
\usepackage{tikz}
\usepackage{pgfplots}
\usepackage{pgfplotstable}
\usetikzlibrary{positioning}
\makeatletter
\definecolor{gnuplot@orange}{RGB}{229,158,0}
\definecolor{gnuplot@purple}{RGB}{148,0,212}
\definecolor{gnuplot@red}{RGB}{200,0,0}
\definecolor{gnuplot@lightblue}{RGB}{87,181,232}
\definecolor{gnuplot@green}{RGB}{0,158,65}
\definecolor{gnuplot@darkblue}{RGB}{0,115,179}
\definecolor{gnuplot@yellow}{RGB}{240,227,66}
\pgfplotscreateplotcyclelist{colorGPL}{%
gnuplot@darkblue,every mark/.append style={fill=gnuplot@darkblue!80!black},mark=o\\%
gnuplot@red,every mark/.append style={fill=gnuplot@red!50!black},mark=square\\%
gnuplot@green,every mark/.append style={fill=gnuplot@green!50!black},mark=triangle*\\%
gnuplot@orange,every mark/.append style={fill=gnuplot@orange!80!black,mark size=2.5pt},mark=x\\%
gnuplot@purple,mark=diamond\\%
black,densely dashed,every mark/.append style={solid,fill=gnuplot@darkblue!80!black},mark=square*\\%
gnuplot@lightblue,densely dashed,every mark/.append style={solid,fill=gnuplot@lightblue!30!black},mark=otimes*\\%
red!80!white,densely dashed,every mark/.append style={solid,fill=red!40!white},mark=oplus*\\%
}
\makeatother

\pgfplotsset{compat=1.9}

\ifpdf
  \DeclareGraphicsExtensions{.eps,.pdf,.png,.jpg}
\else
  \DeclareGraphicsExtensions{.eps}
\fi

\newcommand{\TheTitle}{A Hermite-like basis for faster matrix-free evaluation of interior penalty discontinuous Galerkin operators}

\title{{\TheTitle}\thanks{This
    work was supported by the German Research Foundation (DFG) under
    the project ``High-order discontinuous Galerkin for the exa-scale''
    (ExaDG) within the priority program ``Software for Exascale Computing''
    (SPPEXA), grant agreement no.~KO5206/1-1 and KR4661/2-1, as well as the
Bayerische Kompetenznetzwerk f\"ur Technisch-Wissenschaftliches Hoch- und
H\"ochstleistungsrechnen (KONWIHR) in the framework of the project
``Performance tuning of high-order discontinuous Galerkin solvers for
SuperMUC-NG''. The authors gratefully acknowledge the Gauss Centre for Supercomputing e.V.~(\texttt{www.gauss-centre.eu}) for funding this project
by providing computing time on the GCS Supercomputer SuperMUC at Leibniz Supercomputing Centre (LRZ, \texttt{www.lrz.de})
through project id pr83te.}}

\author{ Martin Kronbichler\thanks{Institute for Computational Mechanics,
    Technical University of Munich, Boltzmannstr.~15, 85748 Garching, Germany (\texttt{\{kronbichler,fehn,munch\}@lnm.mw.tum.de}).}%
  \and%
  Katharina Kormann\thanks{Max Planck Institute for Plasma Physics,
    Boltzmannstr.~2, 85748 Garching, Germany
    (\texttt{katharina.kormann@ipp.mpg.de}).}%
  \and%
  Niklas Fehn\footnotemark[2]%
  \and%
  Peter~Munch\footnotemark[2]~\thanks{Institute of Materials Research, Materials Mechanics, Helmholtz-Zentrum Geesthacht, Max-Planck-Str.~1, 21502 Geesthacht, Germany (\texttt{peter.muench@hzg.de}).}\and%
  Julius Witte\thanks{Interdisciplinary Center for Scientific Computing,
    Heidelberg University, Im Neuenheimer Feld 205, 69120 Heidelberg,
    Germany (\texttt{julius.witte@iwr.uni-heidelberg.de}).}%
}

\usepackage{amsopn}

\begin{document}

\maketitle

\begin{abstract}
This work proposes a basis for improved throughput of matrix-free evaluation of
discontinuous Galerkin symmetric interior penalty discretizations on hexahedral
elements. The basis relies on ideas of Hermite polynomials.
It is used in a fully discontinuous setting not for higher order continuity but to minimize
 the effective stencil width, namely to limit the
neighbor access of an element to one data point for the function value and one
for the derivative. The basis is extended to higher orders with nodal
contributions derived from roots of Jacobi polynomials and extended to multiple
dimensions with tensor products, which enable the use of sum factorization.
The beneficial effect of the
reduced data access on modern processors is shown.
Furthermore, the viability of the basis in the context of
multigrid solvers is analyzed. While a plain
point-Jacobi approach is less efficient than with the best nodal polynomials,
a basis change via sum-factorization techniques enables the combination of
the fast matrix-vector products with effective multigrid constituents. The
basis change is essentially for free on modern hardware because these
computations can be hidden behind the cost of the data access.
\end{abstract}

\noindent \textbf{Key words.}
  High-order discontinuous Galerkin, matrix-free method, sum factorization,
  high-performance computing, multigrid method.

\section{Introduction}
\label{sec:intro}

High-order discontinuous Galerkin (DG) methods, traditionally most popular for
conservation laws, are gaining in popularity also for elliptic problems. For
example, representing the pressure Poisson equation in incompressible flows with
discontinuous pressure spaces combines naturally with $H$(div) conforming
velocity spaces \cite{Lehrenfeld16} (e.g., Raviart--Thomas spaces in the
case of hexahedral elements) or stabilized DG spaces \cite{Fehn18a,Piatkowski18}.
In the context of solving linear systems, DG methods are often
considered expensive due to a
relatively wide stencil, because each shape function is densely
coupled to all the shape functions within a cell as well as to the neighbors.
In primal formulations of elliptic operators, such
as the interior penalty method for the scalar Poisson equation considered in this work,
the neighbor coupling links to all shape functions which have a
non-zero value or non-zero first derivative on the face. This corresponds to all unknowns
on the neighboring elements for conventional bases.
For tensor product shape functions, the polynomial space is
also much larger than necessary to describe the complete polynomial space of
degree $p$. In order to reduce the memory requirements and memory transfer of matrix-based approaches, more compact schemes like the compact discontinuous Galerkin
method~\cite{Peraire07}, the hybridizable discontinuous Galerkin method
\cite{Cockburn09}, or the line discontinuous Galerkin method \cite{Persson13}
have been developed. Also, finite element bases with reduced sparsity within
an element have been proposed~\cite{Beuchler12}.

An alternative path to more efficient DG solvers are matrix-free
implementations, using ideas originally developed within spectral elements
\cite{Deville02,Kopriva09,O80}. Rather than assembling a global sparse matrix
that is later used in some iterative linear solver, these methods compute the
matrix-vector product in terms of the underlying integrals on the fly by fast
quadrature. These matrix-free methods allow to use the tensor
product (Kronecker) structure in the shape functions and the quadrature points:
For interpolation
between solution coefficients and values or derivatives at quadrature points,
the so-called sum-factorization approach separates the
interpolation along each of the coordinate directions from the constant factors
in the other directions,
thereby reaching a complexity of
$\mathcal O(p^{d+1})$ arithmetic operations per cell in the polynomial degree
$p$ in $d$ dimensions for the interpolation. Scaled to the number of operations per degree
of freedom, the cost is linear, $\mathcal O(p)$.
The observed throughput of operator evaluation with
sum-factorization schemes is often constant per degree of freedom for moderately
high polynomial degrees $p\leq 10$ when face integrals and memory access of complexity $\mathcal O(p^d)$
dominate over the $\mathcal O(p^{d+1})$ complexity of cell integrals~\cite{Fischer19}.
Even better, when
executed on modern memory-bandwidth-starved hardware, matrix-free operator
evaluation for polynomial degrees $3\leq p\leq 10$ is also several times faster
than a sparse matrix-vector product for \emph{continuous linear} finite elements with the same
number of unknowns \cite{Kronbichler18}.
Furthermore, the final matrix entries do not possess a tensor product
structure for variable coefficients
or deformed elements and thus the sparse matrix-vector product involves $\mathcal O(p^{d})$
operations per unknown.

While fast operator evaluation is well-established in explicit time integration,
iterative solvers are often preconditioned by schemes such as the Gauss--Seidel
relaxation or incomplete factorizations, which explicitly rely on the matrix
entries and are not compatible with integration-based matrix-free methods. For
multigrid schemes, selected smoothers such as a Chebyshev iteration around the
Jacobi method \cite{Adams03,Sundar15} or methods based on the fast
diagonalization method \cite{lottes05} enable optimal $\mathcal O(1)$ storage
complexity per unknown and ensure that the matrix-vector product dominates.
Apart from the serial performance, they
are particularly interesting in a parallel setting.
We note that the application metric driving the present work is time to solution, given by the
number of effective iterations (matrix-vector products) times the cost of a
matrix-vector product as discussed in \cite{Fehn18b}.

Regarding the efficiency of matrix-free methods, DG
spectral-ele\-ment methods (DGSEM) with collocated node points for Lagrange
polynomials and integration points have been most popular~\cite{Kopriva09}.
For DG methods, one can choose between nodal polynomials in the points of Gauss--Lobatto
polynomials and in the points of the Gauss
quadrature. While the latter ensures exact integration on affine geometries, the former has nodes at the boundary which simplifies access to values on faces, at the price of additional discretization errors due to the Gauss--Lobatto quadrature.

Progress in computer hardware has made computations
cheaper relative to data movement. As a consequence, the data access patterns
are as important as the number of arithmetic operations in the selection of
a basis, as has been found in a recent study~\cite{Schoeder18}. This is
because matrix-free operator evaluation with a state-of-the-art implementation, like the one from the deal.II finite element library~\cite{dealII90} used in this work, comes close to the
throughput of simply streaming the input and output vectors on modern processors.
Computer architecture therefore fosters a data-centric point of view
with the goal to reduce the memory transfer to a minimum when traversing
through the grid, an optimization technique well-known from high-performance
finite difference implementations~\cite{Hager11}. The present work meets
these requirements by proposing a Hermite-like basis compatible with a single
sweep through data for both cell and face integrals with minimal data access
for the prototype second-derivative DG operator, the symmetric interior penalty
method. The selection of
this type of shape functions is orthogonal to the typical use of Hermite
polynomials in the context
of continuous finite elements for reaching higher order of regularity (e.g.
$C^1$ continuity), see also a recent work on wave propagation \cite{Appelo18}.
Instead, we keep a fully discontinuous $L_2$-conforming ansatz
space similar to the work by \cite{Hemker03}. The chosen basis covers the same function space as a nodal Lagrange basis
and is integrated with the same Gaussian quadrature producing exact integrals on affine element
shapes, and, therefore, does not alter accuracy.

This work is structured as follows. The DG discretization is introduced in
section~\ref{sec:model} and the chosen matrix-free implementation in
section~\ref{sec:matrixfree}. The Hermite-like basis functions with
well-conditioned interpolation are constructed in
section~\ref{sec:polynomials}.
Section~\ref{sec:experiments} gives an in-depth
performance analysis of the matrix-free operator evaluation with the new basis against established approaches, supported by a cache analysis. In section~\ref{sec:multigrid} an efficient multigrid
scheme with this method is discussed before conclusions are given in section~\ref{sec:conclusions}.

\section{Discontinuous Galerkin discretization of the Laplacian}
\label{sec:model}

We consider
the symmetric interior penalty discretization of the scalar
Laplacian according to \cite{Arnold02}, whose homogeneous part on an element $K$ is given by
the weak form
\begin{equation}\label{eq:discrete_laplace}
\begin{aligned}
  &\left(\nabla \varphi_i, \nabla u_h\right)_{K}
  +\sum_{F\in \text{faces}(K)}\hspace{-0.3cm}
  \left<\varphi_i \hat{n}, \sigma [\![u_h]\!]\right>_{F}
    -\left<\varphi_i \hat{n}, \{\!\{\nabla u_h\}\!\}\right>_{F}
    -\left<\frac{\nabla \varphi_i}{2}, [\![u_h]\!]\right>_{F},
\end{aligned}
\end{equation}
where $u_h = \sum_{j} \varphi_j u_j$ is the finite element interpolation based
on shape functions $\varphi_j$ and the vector entries $u_j$. In
equation~\eqref{eq:discrete_laplace}, the operator
$\{\!\{v\}\!\}=\frac{v^-+v^+}{2}$ denotes the average of the quantity $v$ over
a face from the interior value $v^-$ on the cell $K$ and the value $u^+$ on the
neighbor $K^+$ behind the respective face, and
$[\![v]\!]=\hat n^- v^- + \hat n^+ v^+ = \hat n (v^- - v^+)$ the directed jump over the
interface along the direction of the outer normal $\hat n=\hat n^-$ of cell $K$. The
parameter $\sigma$ is a penalty parameter chosen large enough in terms of the
polynomial degree $p$ and the mesh size $h$ to make the
final weak form coercive, $\sigma = \nicefrac{(p+1)^2}{h}$,
see e.g.~\cite{Epshteyn07}. The bilinear forms $(\cdot, \cdot)_K$ and
$\left<\cdot,\cdot\right>_F$ denote the integration of the product of the two
arguments over the cell $K$ and on the face $F$, respectively. Boundary conditions
are assumed to be implemented via the mirror principle, e.g., $u^+ = -u^-$ on
homogeneous Dirichlet boundaries.

In this work, we consider high-order DG methods on a mesh of hexahedral
elements where sum factorization is most straight-forward.
In reference
coordinates, the basis functions
\begin{equation}\label{eq:phi_tensor}
\varphi_i(\boldsymbol{\xi}) = \phi_{i_1}(\xi_1) \phi_{i_2}(\xi_2)\ldots
\phi_{i_d}(\xi_d)
\end{equation}
are constructed as the tensor product of one-dimensional polynomials $\phi$ in
each of the $d$ coordinate directions. For polynomial degree $p$,
this results in $(p+1)^d$ basis functions per element. In the usual finite
element setting, the functions in real space~$\boldsymbol x$ are defined
by a transformation from a reference space
$\boldsymbol \xi\in [0,1]^d$ to the real space using some polynomial mapping.
In the DG method, no continuity over the element faces is imposed. The terms in
\eqref{eq:discrete_laplace} are tested by all test functions~$\varphi_i$, $i=1,
\ldots, (p+1)^d$, and for all elements $K$ in the mesh,  resulting in a linear system
with $n_\text{DoF} = n_\text{elements} (p+1)^d$ global equations and unknowns.

\section{Matrix-free evaluation of finite element operators}
\label{sec:matrixfree}

In matrix-free operator evaluation, the matrix-vector product
$\boldsymbol{y} = A \boldsymbol u$ is evaluated in terms of the weak
form~\eqref{eq:discrete_laplace} using the finite element field $u_h$
associated with the vector values in $\boldsymbol u$. The integrals are
computed by quadrature on $n_q^d$ Gauss quadrature points.
On affine geometries and
for constant coefficients, choosing $n_q= p+1$ quadrature points ensures exact
integration. For variable coefficients, nonlinear terms, or curved geometries,
$n_q>p+1$ is sometimes needed for accurate results. For the example of
the cell term $(\nabla \varphi_i, \nabla u_h)_K$, the integral is approximated by
\begin{equation}\label{eq:integral_laplace}
  \begin{aligned}
    (\nabla \varphi_i, \nabla u_h)_{K}&=\int_{\Omega_\text{unit}} \left(\mathcal J_{K}^{-\mathsf T} \nabla _{\boldsymbol \xi} \varphi_i \right)^{\mathsf T} \left(\mathcal J_{K}^{-\mathsf T} \sum_{j=1}^{(p+1)^d}\nabla _{\boldsymbol \xi} \varphi_j u_{j}^{(K)}\right)  \mathrm{det}(\mathcal J_{K})\, \mathrm{d}\boldsymbol{\xi}\\
    &\approx\sum_{q=1}^{n_q^d}(\nabla _{\boldsymbol \xi} \varphi_i(\boldsymbol \xi_q))^{\mathsf T} \underbrace{\mathcal J_{K}^{-1}  \mathcal J_{K}^{-\mathsf T} \mathrm{det}(\mathcal J_{K}) w_q}_{\text{depends only on q-point}} \sum_{j=1}^{(p+1)^d}\nabla _{\boldsymbol \xi} \varphi_j(\boldsymbol \xi_q) u_{j}^{(K)}.
  \end{aligned}
\end{equation}
Here, $\mathcal J_{K}$ denotes the Jacobian of the mapping from the reference to the real cell at the quadrature points and $w_q$ the quadrature weight.
Since the metric terms
do not depend on the shape function indices $i$ and~$j$, they
are evaluated outside the~$i$ and~$j$ loops, a common
abstraction in matrix-free methods \cite{Fischer19,Knepley13,Kronbichler12,Kronbichler19}.

In this version of the algorithm, the work at quadrature points
is of complexity $\mathcal O\left(n_q^d\right)$, whereas the interpolation
sum over $j$ as well as testing by all test functions $\varphi_i$ would both imply a cost of
$\mathcal O\left((p+1)^d n_q^d\right)$ per element. The interpolation
is amenable to BLAS-3 linear algebra because the same reference-element
operations is applied on each element $K$. Furthermore, the interpolation and
derivative matrices have more structure in case the
polynomials are the tensor product of 1D polynomials according
to~\eqref{eq:phi_tensor} and integrated with a tensor product quadrature formula:
The two nested loops over $(p+1)^d$ basis functions and $n_q^d$ points can be broken down into a
series of smaller loops along each of the coordinate directions by a technique
called sum factorization~\cite{Deville02,Kopriva09}. The overall complexity per element is
then $\mathcal O(p^{d+1})$ for $n_q\sim p$.

Let us denote by $S=S_{q,j}$ the $n_q\times (p+1)$ matrix of all 1D shape
functions~$\phi_j$ evaluated at all 1D quadrature points $\xi_q$ and by
$D=D_{q,j}$ the derivative of the shape functions
$\phi_{j}^\text{1D,co}(\xi_q)$ with nodes in the quadrature points evaluated at
$\xi_q$ (collocation basis). Then, the sum over $j$
evaluated at all quadrature points in~\eqref{eq:integral_laplace}
can be written in matrix-vector notation as
\begin{equation}\label{eq:tensor_cell}
    \begin{bmatrix}
    \nicefrac{\partial \boldsymbol u_h}{\partial \xi_1}\\
    \nicefrac{\partial \boldsymbol u_h}{\partial \xi_2}\\
    \vdots\\
    	\nicefrac{\partial \boldsymbol u_h}{\partial \xi_d}
    	\end{bmatrix}
    = \begin{bmatrix}
      I \otimes \ldots \otimes I \otimes D\\
      I \otimes \ldots \otimes D \otimes I\\
      \vdots\\
      D \otimes I \otimes \ldots \otimes I
    \end{bmatrix}
    \begin{bmatrix}
      S \otimes \ldots \otimes S \otimes S
    \end{bmatrix}
  \boldsymbol u^{(K)}.
\end{equation}
The first multiplication
$\begin{bmatrix} S \otimes S \otimes \ldots \otimes S
\end{bmatrix} \boldsymbol u^{(K)}$ is a basis change, going from the values~$\boldsymbol u^{(K)}$ in the solution vector to a representation in the
nodal basis of Lagrange polynomials defined in the quadrature points. The
gradient operation is then performed in this basis.
In equation~\eqref{eq:tensor_cell}, the evaluation of the
matrices in Kronecker product form is understood by the usual 1D contractions
along each dimension with sum factorization, which amounts to $2d$ products for
the complete gradient~\eqref{eq:tensor_cell}, see~\cite{Fischer19,Kronbichler19}
and references therein. The
integration step involves the transpose of the matrix from~\eqref{eq:tensor_cell}.

Analogously, the evaluation of the face integrals of
equation~\eqref{eq:discrete_laplace} involves reference-cell operations from both
sides of a face, as well as geometrical and equation-dependent operations at
quadrature points. For a face with normal vector in negative $\xi_1$ direction
in three dimensions, the interpolation is given by
\begin{equation}\label{eq:tensor_face}
    \begin{bmatrix}
    \boldsymbol u_h\\
    \nicefrac{\partial \boldsymbol u_h}{\partial \xi_2}\\
    \nicefrac{\partial \boldsymbol u_h}{\partial \xi_3}\\
    	\nicefrac{\partial \boldsymbol u_h}{\partial \xi_1}
    	\end{bmatrix}
    	=
    \begin{bmatrix}
      \begin{bmatrix}
    \begin{bmatrix}
      I \otimes I\\
      I \otimes D\\
      D \otimes I
    \end{bmatrix}
        \begin{bmatrix}
      S \otimes  S
    \end{bmatrix}\end{bmatrix}
    &  0 \\
    0 & \begin{bmatrix}S\otimes S \end{bmatrix}
  \end{bmatrix}
    \begin{bmatrix}
       I \otimes I \otimes S_{\mathsf{f}}  \\
       I \otimes I \otimes D_{\mathsf{f}}
    \end{bmatrix}
  \boldsymbol u^{(K)},
\end{equation}
where the $(p+1)\times 1$ matrix $S_\mathsf{f}$ contains the values of all shape
functions $\phi_i$ at $\xi=0$, and the $(p+1)\times 1$ matrix $D_\mathsf{f}$
contains the evaluation of the first derivative of all shape functions $\phi_i'$
at $\xi=0$. Note that this operation returns four vectors of length~$n_q^2$,
representing the values $u_h(\boldsymbol \xi_q)$ and the three components of the
gradient $\nabla_{\boldsymbol \xi} u_h(\boldsymbol \xi_q)$ in reference coordinates. Similar operations
can be defined for the other faces and, in transposed form, for integration.
The key observation for face integrals is that the initial multiplications by the
matrices $S_{\mathsf{f}}$ and $D_{\mathsf{f}}$ extract the degrees of freedom
relevant for face integrals from
the solution vector $\boldsymbol u^{(K)}$. Entries where both $S_{\mathsf{f}}$
and~$D_{\mathsf{f}}$ are zero need not be loaded, independent of the geometry.
As proposed in~\cite{Kronbichler19}, interpolation operations on interior faces
can share some basis change operations with the cell integrals,
such that the full interpolation according to equation
\eqref{eq:tensor_face} is only necessary on the neighbors' face data $u^+_h$.

In terms of the arrangement of cell and face integrals, we concentrate on an
algorithm that is referred to as an ``element-wise evaluation of face integrals''
in \cite{Kronbichler19}. Here, all integrals related to an element $K$ (or to a
batch of elements when vectorizing over several elements) are computed together.
This enables a single write operation to the result vector with a single sweep through
the vector and thus minimizes
the data movement. At the same time, this setup evaluates the flux on interior faces
twice. For a comparison to other arrangements of face integrals and their properties, we refer to~\cite[Sec.~4.6]{Kronbichler19}.

In this work, we assume that the geometry terms $\mathcal J$ are possibly
space-varying, necessitating separate interpolation and integration steps
with all cross-terms in the derivatives.
If $\mathcal J$ were constant within an element, i.e., a cell subject to an affine
transformation,
the integration could be completed in reference coordinates and the full
interpolation and integration step could be expressed by a sum of a few
Kronecker matrices, reducing the number of arithmetic operations especially on faces.
We do not consider those optimizations in this work because pre-computed Kronecker
matrices would allow for further algorithmic rearrangements such as element-wise
static condensation \cite{Huismann17}, and are not applicable to variable
geometries for which DG methods have been originally developed.

\section{Hermite-like basis functions}
\label{sec:polynomials}

The selection of the basis in DG methods is more relaxed as compared with
continuous finite elements where inter-element continuity is imposed by
shared nodes. The following three principles can be applied:
\begin{itemize}
\item a modal basis where the final mass matrix is diagonal to ensure cheap
explicit time integration, see e.g.~\cite{Crivellini11}, also easily separating low/high frequency content,
\item a basis where the transformation $S$ is the identity matrix, i.e., a nodal
basis with the nodes coinciding with the points of the quadrature formula, the
typical spectral element/DGSEM setup \cite{Hindenlang12,Kopriva09}, or
\item a basis where the number of entries in $S_\mathsf{f}$ and $D_\mathsf{f}$
is explicitly minimized.
\end{itemize}
The collocation setup with $S=I$ and $n_q = p+1$ Gaussian quadrature points is
widely used and
simplifies several operations in both the cell integrals of
equation~\eqref{eq:tensor_cell}
as well as face integrals of equation~\eqref{eq:tensor_face}. Furthermore, it
also implies a diagonal mass matrix.
However, for Gaussian quadrature all entries in
$S_\mathsf{f}$ and $D_\mathsf{f}$ are non-zero because all roots of the underlying
Legendre polynomials are strictly inside the 1D reference element. When switching
to the
less accurate collocated Gauss--Lobatto quadrature, $S_\mathsf{f}=[1,0,\ldots,0]$ becomes
simple, but $D_\mathsf{f}$ is still dense. Furthermore, as soon as
over-integration with $n_q > p+1$ is used, the identity $S=I$ is lost and all
matrices in the formulas \eqref{eq:tensor_cell} and
\eqref{eq:tensor_face} must be used.

\subsection{Construction of the basis}
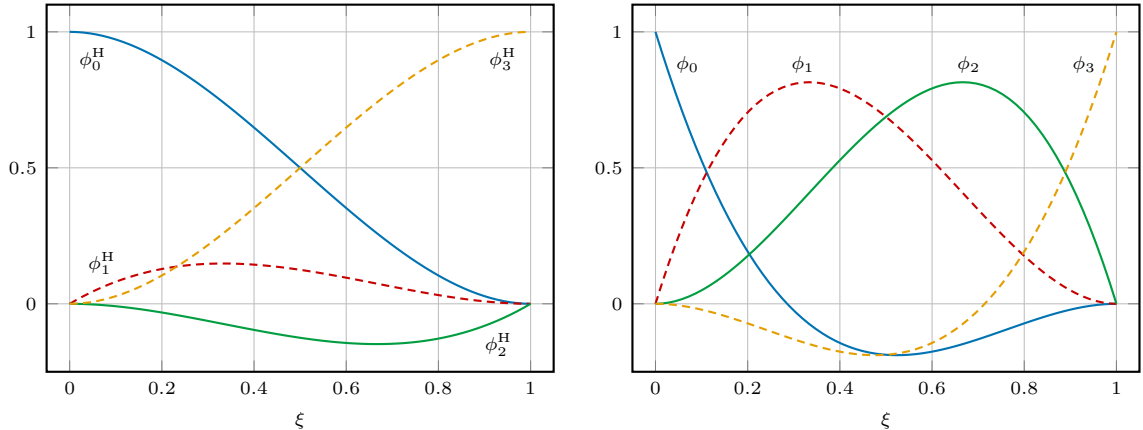
\begin{figure}
  \centering
  \begin{tikzpicture}[domain=0:1]
    \begin{axis}[
      width=0.55\textwidth,
      height=0.43\textwidth,
      xlabel={$\xi$},
      title style={font=\scriptsize},
      tick label style={font=\scriptsize},
      label style={font=\scriptsize},
      legend style={font=\scriptsize},
      xmin=-0.05,xmax=1.05,
      ymin=-0.25,ymax=1.1,
      grid,
      cycle list name=colorGPL,
      mark size=1.8,
      thick
      ]
      \addplot[gnuplot@darkblue,samples=70] (\x, {(2*\x+1)*(\x-1)*(\x-1)});
      \node at (axis cs:0.05,0.9) {\scriptsize $\phi_0^\text{H}$};
      \addplot[gnuplot@red,densely dashed,samples=70] (\x, {(\x)*(\x-1)*(\x-1)});
      \node at (axis cs:0.07,0.15) {\scriptsize $\phi_1^\text{H}$};
      \addplot[gnuplot@green,samples=70] (\x, {(\x)*(\x)*(\x-1)});
      \node at (axis cs:0.93,-0.15) {\scriptsize $\phi_2^\text{H}$};
      \addplot[gnuplot@orange,densely dashed,samples=70] (\x, {-(2*\x-3)*(\x)*(\x)});
      \node at (axis cs:0.94,0.9) {\scriptsize $\phi_3^\text{H}$};
    \end{axis}
  \end{tikzpicture}
  \hfill
  \begin{tikzpicture}[domain=0:1]
    \begin{axis}[
      width=0.55\textwidth,
      height=0.43\textwidth,
      xlabel={$\xi$},
      title style={font=\scriptsize},
      tick label style={font=\scriptsize},
      label style={font=\scriptsize},
      legend style={font=\scriptsize},
      xmin=-0.05,xmax=1.05,
      ymin=-0.25,ymax=1.1,
      grid,
      cycle list name=colorGPL,
      mark size=1.8,
      thick
      ]
      \addplot[gnuplot@darkblue,samples=70] (\x, {-3.5*(\x-2/7)*(\x-1)*(\x-1)});
      \node at (axis cs:0.07,0.88) {\scriptsize $\phi_0$};
      \addplot[gnuplot@red,densely dashed,samples=70] (\x, {5.5*(\x)*(\x-1)*(\x-1)});
      \node at (axis cs:0.32,0.88) {\scriptsize $\phi_1$};
      \addplot[gnuplot@green,samples=70] (\x, {-5.5*(\x)*(\x)*(\x-1)});
      \node at (axis cs:0.68,0.88) {\scriptsize $\phi_2$};
      \addplot[gnuplot@orange,densely dashed,samples=70] (\x, {3.5*(\x-5/7)*(\x)*(\x)});
      \node at (axis cs:0.93,0.88) {\scriptsize $\phi_3$};
    \end{axis}
  \end{tikzpicture}

  \caption{Classical Hermite basis functions \eqref{eq:hermite_basic} (left) and proposed Hermite-like basis functions (right) for degree $p=3$.}
  \label{fig:hermite_3}
\end{figure}

In this work, we minimize the access into~$\boldsymbol u$ for face
integrals, i.e., maximize the number of columns where both the matrices $S_\mathsf{f}$ and
$D_\mathsf{f}$ have zero entries.
More precisely, the idea is to generalize the favorable access pattern of Lagrange polynomials in the Gauss--Lobatto points for hyperbolic
problems, where solution values on faces access only $(p+1)^{d-1}$ entries due to $S_\mathsf{f}=[1,0,\ldots,0]$, to settings where values and gradients are needed on faces.
In the recent analysis
\cite{Schoeder18} it has been shown that a tight data access via $S_\mathsf{f}=[1,0,\ldots,0]$ outweighs the
additional calculations with dense matrices $S$ on modern hardware.
With respect to minimizing the number of nonzeros in
$D_\mathsf{f}$, it is necessary to make the first derivative vanish at $\xi=0$
and $\xi=1$ for all shape functions but one, naturally leading to the Hermite
basis for $p=3$ with
\begin{equation}\label{eq:hermite_prop}
D_\mathsf{f}=[0, 1, 0, \ldots, 0].
\end{equation}
The four cubic Hermite polynomials---shown in Figure~\ref{fig:hermite_3} (left)---are defined as
\begin{equation}\label{eq:hermite_basic}
\phi_0^\text{H}(\xi)=2\xi^3-3\xi^2+1,\ \phi_1^\text{H}(\xi) = \xi^3-2\xi^2+\xi,\
\phi_2^\text{H}(\xi) =\xi^3-\xi^2 ,\ \phi_3^\text{H}(\xi)=-2\xi^3+3\xi^2.
\end{equation}
Regarding the case
$p>3$, the support of higher derivatives on the reference cell boundary can
be limited, which leads to the Bernstein--B{\'e}zier polynomial basis. While
this basis has attractive properties such as a fast recursive evaluation of
polynomials or applicability of sum factorization algorithms also for triangular
and tetrahedral elements~\cite{Ainsworth11} including inverse mass
matrices~\cite{Kirby17}, the exponential increase in the condition number requires
careful algorithm selection to not spoil the attractive properties by roundoff
effects. Hermite-type polynomials have already been considered in the DG
context for constructing better multigrid smoothers in~\cite{Hemker03}.

If we limit ourselves to minimizing the entries of $S_\mathsf{f}$ and $D_\mathsf{f}$,
higher order polynomials can be constructed via
\begin{equation}\label{eq:hermite_legendre}
\phi_{i}^{(p)}(\xi) = 16\xi^2 (1-\xi)^2 \sigma_{i}^{(p)}(\xi),
\end{equation}
where $\sigma_{i}^{(p)}(\xi)$ with $i=4, \ldots, p$ is a polynomial of degree up
to $p-4$. These polynomials are bubble functions with vanishing values and first derivatives
on both faces. A naive possibility is to choose Legendre polynomials.
Table~\ref{tab:condition_mass} lists the quality of this
basis measured by the condition number of the 1D consistent mass matrix.
The condition number is $10^3$ already for the cubic Hermite case $p=3$
and further deteriorates as the degree increases.

\begin{table}
  \caption{Condition number of the 1D mass matrix for the Hermite interpolation
extended by Legendre polynomials in~\eqref{eq:hermite_legendre} and for the
proposed Hermite-like interpolation. The condition number of the standard nodal
Lagrange basis on the Gauss--Lobatto points is included for reference.}
\label{tab:condition_mass}
\footnotesize
\centering
\begin{tabular}{lccc}
\hline
$p$ & Hermite + Legendre extension & Hermite-like basis & Lagrange basis on GL points \\
\hline
$3$ & $1.06\cdot 10^3$ & $17.2$ & $8.65$\\
$4$ & $6.58\cdot 10^3$ & $16.8$ & $10.4$\\
$5$ & $1.88\cdot 10^4$ & $16.0$ & $12.0$\\
$6$ & $6.03\cdot 10^4$ & $16.3$ & $13.6$\\
$7$ & $1.28\cdot 10^5$ & $17.1$ & $15.2$\\
$8$ & $2.85\cdot 10^5$ & $18.2$ & $16.7$\\
$10$ & $9.76\cdot 10^5$ & $20.7$ & $19.9$\\
$15$ & $9.43\cdot 10^6$ & $27.9$ & $27.7$\\
$20$ & $5.67\cdot 10^7$ & $35.7$ & $35.5$\\
$25$ & $2.22\cdot 10^8$ & $43.5$ & $43.4$\\
$30$ & $7.73\cdot 10^8$ & $51.6$ & $51.2$\\
\hline
\end{tabular}
\end{table}

\begin{figure}
\strut\hfill
  \begin{tikzpicture}[domain=0:1]
   \begin{axis}[
     width=0.33\textwidth,
     height=0.4\textwidth,
     xlabel={$\xi$},
     title style={font=\scriptsize},
     tick label style={font=\scriptsize},
     label style={font=\scriptsize},
     legend style={font=\scriptsize},
     xmin=-0.05,xmax=1.05,
     ymin=-0.38,ymax=1.1,
     xtick={0,0.2,0.4,0.6,0.8,1},
     grid,
     cycle list name=colorGPL,
     mark size=1.8,
     thick
     ]
     \addplot[gnuplot@darkblue,samples=70] (\x, {(1-\x)});
      \node at (axis cs:0.28,0.85) {\scriptsize $\phi_0$};
     \addplot[gnuplot@orange,densely dashed,samples=70] (\x, {\x});
      \node at (axis cs:0.72,0.85) {\scriptsize $\phi_1$};
   \end{axis}
  \end{tikzpicture}
  \hfill
  \begin{tikzpicture}[domain=0:1]
   \begin{axis}[
     width=0.33\textwidth,
     height=0.4\textwidth,
     xlabel={$\xi$},
     title style={font=\scriptsize},
     tick label style={font=\scriptsize},
     label style={font=\scriptsize},
     legend style={font=\scriptsize},
     xmin=-0.05,xmax=1.05,
     ymin=-0.38,ymax=1.1,
     xtick={0,0.2,0.4,0.6,0.8,1},
     grid,
     cycle list name=colorGPL,
     mark size=1.8,
     thick
     ]
     \addplot[gnuplot@darkblue,samples=70] (\x, {(\x-1)*(\x-1)});
      \node at (axis cs:0.14,0.92) {\scriptsize $\phi_0$};
     \addplot[gnuplot@red,densely dashed,samples=70] (\x, {-2*(\x)*(\x-1)});
      \node at (axis cs:0.5,0.58) {\scriptsize $\phi_1$};
     \addplot[gnuplot@orange,densely dashed,samples=70] (\x, {\x*\x});
      \node at (axis cs:0.86,0.92) {\scriptsize $\phi_2$};
   \end{axis}
  \end{tikzpicture}
  \hfill
  \begin{tikzpicture}[domain=0:1]
    \begin{axis}[
      width=0.49\textwidth,
      height=0.4\textwidth,
      xlabel={$\xi$},
      title style={font=\scriptsize},
      tick label style={font=\scriptsize},
      label style={font=\scriptsize},
      legend style={font=\scriptsize},
      xmin=-0.05,xmax=1.05,
      ymin=-0.38,ymax=1.1,
      xtick={0,0.2,0.4,0.6,0.8,1},
      grid,
      cycle list name=colorGPL,
      mark size=1.8,
      thick
      ]
      \addplot[gnuplot@darkblue,samples=150] (\x,
      {-429*(\x-0.0606)*(\x-0.201)*(\x-0.396)*(\x-0.604)*(\x-0.799)*(\x-1)*(\x-1)}); 
      \addplot[gnuplot@red,densely dashed,samples=150] (\x, {751*(\x)*(\x-0.201)*(\x-0.396)*(\x-0.604)*(\x-0.799)*(\x-1)*(\x-1)});
      \addplot[gnuplot@yellow,samples=150] (\x, {-825*(\x)*(\x)*(\x-0.396)*(\x-0.604)*(\x-0.799)*(\x-1)*(\x-1)});
      \addplot[black,densely dashed,samples=150] (\x, {1073*(\x)*(\x-0.201)*(\x)*(\x-0.604)*(\x-0.799)*(\x-1)*(\x-1)});
      \addplot[gnuplot@lightblue,samples=150] (\x, {-1073*(\x)*(\x-0.201)*(\x-0.396)*(\x)*(\x-0.799)*(\x-1)*(\x-1)});
      \addplot[gnuplot@purple,densely dashed,samples=150] (\x, {825*(\x)*(\x-0.201)*(\x-0.396)*(\x-0.604)*(\x)*(\x-1)*(\x-1)});
      \addplot[gnuplot@green,samples=150] (\x, {-751*(\x)*(\x-0.201)*(\x-0.396)*(\x-0.604)*(\x-0.799)*(\x)*(\x-1)});
      \addplot[gnuplot@orange,densely dashed,samples=150] (\x, {429*(\x)*(\x-0.201)*(\x-0.396)*(\x-0.604)*(\x-0.799)*(\x)*(\x-0.939)});
    \end{axis}
  \end{tikzpicture}
  \hfill\strut
  \caption{Hermite-like basis for $p=1$ (left), $p=2$ (middle), and $p=7$ (right).}
  \label{fig:hermite_7}
\end{figure}
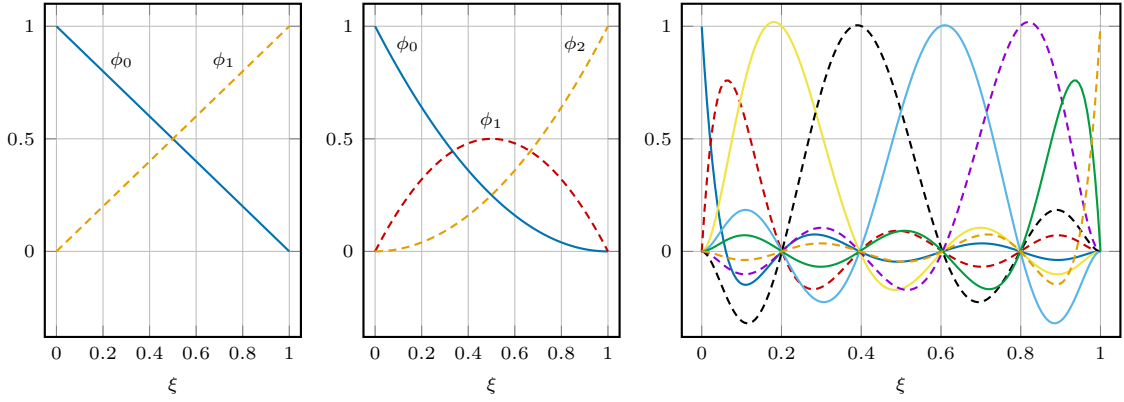

In order to improve the conditioning of interpolation, we relax the Hermite
polynomials into a basis we call ``Hermite-like'', illustrated in Figures~\ref{fig:hermite_3} and \ref{fig:hermite_7}. The construction of the basis is as follows:
\begin{enumerate}
\item \textbf{Relaxation to two nonzero entries in  $D_\mathsf{f}$ for improved conditioning}. In order to improve conditioning with $p=3$, we relax the polynomial
$\phi_0^\text{H}$ by allowing $\phi_0'(0)\neq 0$ via a free root $\xi_1$ in the
form $\phi_0(\xi)
= \alpha_0 (\xi - \xi_1) (\xi-1)^2$ with some constant $\alpha_0$. Furthermore,
we set $\phi_1(\xi) = \alpha_1\xi (\xi-1)^2$ with some constant $\alpha_1$. The
root $\xi_1$ is fixed by the heuristic argument that the orthogonality between
$\phi_0$ and $\phi_1$, $\int_{0}^1 \phi_0 \phi_1 d\xi = 0$, will keep the
condition number low. For $p=3$, we obtain $\xi_1 =
\frac{2}{7}$. The requirement $\phi_0(0)=1$ then gives $\alpha_0 =
-\frac{1}{\xi_1} = -\frac 72$. The second constant $\alpha_1$ is determined by
the condition $\phi_1'(0) = -\phi_0'(0)$, which gives a relaxed condition
\begin{equation}\label{eq:hermite_relaxed}
D_\mathsf{f}=[-\alpha_1, \alpha_1, 0, \ldots, 0].
\end{equation}
The constant evaluates to $\alpha_1 = \frac{2 \xi_1+1}{\xi_1} =
\frac{11}{2}$. The polynomials $\phi_2$ and $\phi_3$ are obtained by mirroring
$\phi_1$ and $\phi_0$ at $\xi = \frac 12$.

Note that relaxing $D_\mathsf{f}$ to two non-zero entries does not worsen data
access because the first entry of $S_\mathsf{f}$ is nonzero anyway and equations
which only need the derivative on faces are uncommon. More
importantly, any non-Cartesian element shape must evaluate tangential
derivatives in reference space as well, for which the values on the face are
needed. For evaluating only the
values on a face, e.g.~for a hyperbolic term,
a single value is touched via~$S_\mathrm{f}$.

\item \textbf{Construction for $p\leq 2$}. We require
$D^\mathsf{f} = [-\alpha_1,\alpha_1,0]$ at $\xi=0$ and
$D^\mathsf{f} = [0,-\alpha_1,\alpha_1]$ at $\xi=1$ with
$\alpha_1=2$ for $p=2$, see Figure~\ref{fig:hermite_7} and
Appendix~\ref{sec:appendix}. For $p=1$,
the standard Lagrange basis $\{1-\xi, \xi\}$ satisfies $S^\mathsf{f} = [1,0]$
and~\eqref{eq:hermite_relaxed}.

\item \textbf{Extension for $p>3$}.
\begin{enumerate}
\item
Higher order polynomials are defined in a nodal way, placing
additional nodes within the reference element. To ensure robust conditioning for
large $p$, the nodes are specified as the roots of the Jacobi
polynomial~$P^{4,4}_{p-3}$, i.e., the Jacobi polynomials orthogonal with respect to the
function $\xi^4(1-\xi)^4$, see also~\cite{Hemker03}. This ensures $L_2$ orthogonality
between the bubble functions \eqref{eq:hermite_legendre}.
Note that we define the reference interval to be $(0,1)$
in this work.
\item The four polynomials active at the boundaries are of degree $p$ with
factors $\xi-\xi_l$ involving the $p-3$ roots of $P^{4,4}_{p-3}$, plus the
boundary contribution $\alpha_0 (\xi - \xi_1) (\xi-1)^2$ for $\phi_0$, for
instance. The additional root
$\xi_1$ is determined by the orthogonality between $\phi_0$ and $\phi_1$. Since
the two polynomials $\phi_0$ and $\phi_1$ only differ by a constant and the
factor $\xi-\xi_1$ or~$\xi$, respectively, the equation for $\xi_1$ is linear
that is easily solved.
\item The weights $\alpha_0$ and $\alpha_1$ of $\phi_0, \phi_1$ are determined
by $\phi_0(0) = 1$, i.e., $S_\mathsf{f} = [1,0,\ldots,0]$, and
$\phi'_0(0) = -\phi'_1(0)$, i.e., relation \eqref{eq:hermite_relaxed}.
This makes sure that the function of all ones, $\boldsymbol
u = [1,1,\ldots,1]^\mathsf T$, represents the constant function $u_h\equiv 1$.
\end{enumerate}
\end{enumerate}

The condition number of the mass matrix resulting from the Hermite-like basis in
Table \ref{tab:condition_mass} shows a
linear increase with the polynomial degree, in analogy to the condition number for the
standard Lagrange basis defined in the points of the Gauss--Lobatto quadrature
formula.
The proposed
construction makes the polynomial basis symmetric with respect to the center of
the reference element $\xi=\frac 12$ in the sense $\phi_i(\xi) = \phi_{p-i}(1-\xi)$ for $i=0,\ldots,p$
like a Lagrange basis. This enables to straight-forwardly apply
optimizations for nodal bases such as the even-odd
decomposition that cuts the work of the 1D operations in
sum factorization into half~\cite{Kronbichler19,Solomonoff92}.

\begin{figure}
\centering
    \begin{tikzpicture}[scale=0.8]
        \foreach \x in{-3, -1.8, -0.6, 0.6}
        \foreach \y in{-5.5, -6.7, -7.9}
        {
          \fill[black!10] (\x,\y) rectangle (\x+1.1,\y+1.1);
          \foreach \xx in {0.05,0.2,0.4,0.7,0.9,1.05}
          \foreach \yy in {0.05,0.2,0.4,0.7,0.9,1.05}
          \draw[black!50] (\x+\xx,\y+\yy) circle (0.05);
        }
        \fill[gnuplot@lightblue!60] (-1.8,-6.7) rectangle(-0.7,-5.6);
        \foreach \x in{-2.95, -2.8, -2.6, -2.3, -2.1, -1.95, -0.55, -0.4, -0.2, 0.1, 0.3, 0.45}
        \foreach \yy in {0.05,0.2,0.4,0.7,0.9,1.05}
        \draw[black,thick] (\x,-6.7+\yy) circle (0.05);
        \foreach \y in{-4.45, -4.6, -4.8, -5.1, -5.3, -5.45, -6.85, -7, -7.2, -7.5, -7.7, -7.85}
        \foreach \xx in {0.05,0.2,0.4,0.7,0.9,1.05}
        \draw[black,thick] (-1.8+\xx,\y) circle (0.05);
        \foreach \xx in {0.05,0.2,0.4,0.7,0.9,1.05}
        \foreach \yy in {0.05,0.2,0.4,0.7,0.9,1.05}
        \draw[black,thick,fill=black] (-1.8+\xx,-6.7+\yy) circle (0.05);
        \node at (-0.3,-8.1)[text width=4.5cm,font=\scriptsize,anchor=north] {2D Lagrange basis in GL points};
        \begin{scope}[xshift=5.5cm]
        \foreach \x in{-3, -1.8, -0.6, 0.6}
        \foreach \y in{-5.5, -6.7, -7.9}
        {
          \fill[black!10] (\x,\y) rectangle (\x+1.1,\y+1.1);
          \foreach \xx in {0.05,0.2,0.4,0.7,0.9,1.05}
          \foreach \yy in {0.05,0.2,0.4,0.7,0.9,1.05}
          \draw[black!50] (\x+\xx,\y+\yy) circle (0.05);
        }
        \fill[gnuplot@lightblue!60] (-1.8,-6.7) rectangle(-0.7,-5.6);
        \foreach \x in{-2.1, -1.95, -0.55, -0.4}
        \foreach \yy in {0.05,0.2,0.4,0.7,0.9,1.05}
        \draw[black,thick] (\x,-6.7+\yy) circle (0.05);
        \foreach \y in{-5.3, -5.45, -6.85, -7.0}
        \foreach \xx in {0.05,0.2,0.4,0.7,0.9,1.05}
        \draw[black,thick] (-1.8+\xx,\y) circle (0.05);
        \foreach \xx in {0.05,0.2,0.4,0.7,0.9,1.05}
        \foreach \yy in {0.05,0.2,0.4,0.7,0.9,1.05}
        \draw[black,thick,fill=black] (-1.8+\xx,-6.7+\yy) circle (0.05);
        \node at (-0.5,-8.1)[text width=4cm,font=\scriptsize,anchor=north] {2D Hermite-like basis functions};
        \end{scope}
        \begin{scope}[xshift=9.2cm,yshift=-6.9cm]
        \filldraw[gnuplot@lightblue!60] (1.3,1.3,2.6) -- (1.3,2.6,2.6) -- (1.3,2.6,1.3) -- (2.6,2.6,1.3) -- (2.6,1.3,1.3)-- (2.6,1.3,2.6) -- (1.3,1.3,2.6);
\foreach \x in{0,1.3,2.6}
{   \draw (0,\x ,2.6) -- (2.6,\x ,2.6);
    \draw (\x ,0,2.6) -- (\x ,2.6,2.6);
    \draw (2.6,\x ,2.6) -- (2.6,\x ,0);
    \draw (\x ,2.6,2.6) -- (\x ,2.6,0);
    \draw (2.6,0,\x ) -- (2.6,2.6,\x );
    \draw (0,2.6,\x ) -- (2.6,2.6,\x );
}
\foreach \x in{1.36,1.54,1.77,2.13,2.36,2.54}
{
\foreach \y in{1.36,1.54,1.77,2.13,2.36,2.54}
{
\draw[black,semithick,fill] (\x,\y,2.6)circle (0.04);
\draw[black,semithick,fill] (2.6,\x,\y)circle (0.04);
\draw[black,semithick,fill] (\x,2.6,\y)circle (0.04);
}
\foreach \y in{1.06,1.24}
{
\draw[black,very thin] (\x,\y,2.6)circle (0.04);
\draw[black,very thin] (\y,\x,2.6)circle (0.04);
\draw[black,very thin] (2.6,\x,\y)circle (0.04);
\draw[black,very thin] (2.6,\y,\x)circle (0.04);
\draw[black,very thin] (\x,2.6,\y)circle (0.04);
\draw[black,very thin] (\y,2.6,\x)circle (0.04);
}
}
        \node at (1.1,-1.23)[text width=3.9cm,font=\scriptsize,anchor=north] {3D Hermite-like basis functions};
        \end{scope}
\end{tikzpicture}
\caption{Data access pattern for $p=5$ in terms of values read (black circles) for the element shaded in blue and values read and written (black disks). Note that the Hermite-like basis is not nodal and the two layers of nodes closest to the surface are highlighted for illustration purposes.}
\label{fig:data_access}
\end{figure}
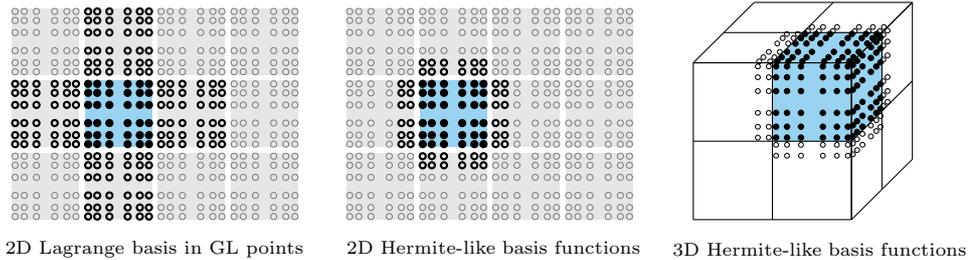

\subsection{Data access}
The data access when computing all cell and face integrals according to \eqref{eq:discrete_laplace} with a matrix-free algorithm is shown in
Figure~\ref{fig:data_access}. From the
figure, it is apparent that the proposed basis minimizes the width of the effective
``element-stencil'', only using one degree of freedom for the values and one for
derivatives in face-normal direction, respectively. Accumulating the access of
an element, the
proposed Hermite-like basis reads
$[(p+1)^d + 4d(p+1)^{d-1}]$ words, compared
with $(2d+1)(p+1)^{d}$ reads for a conventional nodal basis. On general
(curvilinear) element shapes, some geometric information must also be loaded. Computing
a tri-linear mapping adds only little memory transfer for $p\geq
3$ and is of high arithmetic intensity. Isoparametric representations access
$d(p+1)^d$ words per element and significantly increase both the data transfer
and the arithmetic work.\footnote{We note that the Hermite-like basis could also be used to reduce
the data access when evaluating the metric terms $\mathcal J_K$ on neighbors.} To avoid
computations, a separate Jacobian matrix $\mathcal J_K^{-\mathsf T}$ can be stored for
each quadrature point of the cell and $\hat n^\mathsf T \mathcal{J}_K^{-\mathsf{T}}$
for each point on both sides of the faces, for a total of $d^2 (p+1)^d + 4d^2(p+1)^{d-1}$
words per element \cite{Kronbichler19}. Ignoring possible caching for
neighboring vector entries, the proposed
Hermite-like basis improves data access by a ratio 1:7 ($=1+2d$) in the limit of large $p$
when done on a 3D affine mesh without
big geometry data and by a ratio 10:16 on a deformed mesh with separate Jacobians for
each quadrature point, or somewhere in between depending on the evaluation of
the geometry. As some neighboring data can be cached, the difference is smaller
in practice when run on CPUs. As will be shown in the next section, the
difference is significant nonetheless. For a GPU implementation, usually no
neighbor data is cached and the improvement would be more significant.

\section{Efficiency of matrix-free operator evaluation}
\label{sec:experiments}

In order to characterize the matrix-free operator evaluation with
the proposed basis, we consider a 3D benchmark test of the
Laplace operator with homogeneous Dirichlet boundary conditions. We
construct an affine geometry by deforming the brick
$(-0.95,0.95)\times(-0.9,0.89)\times(-0.85,0.83)$ through a linear
transformation with the Jacobian
\[
\mathcal J = \begin{pmatrix}1.12 & 0.24 & 0.36\\0.24 &1.36 &0.48\\ 0.36 & 0.48 &
1.60\end{pmatrix}.
\]
As mentioned above, we perform
integration as if the geometry were deformed,
but using the same merged tensor $\mathcal J_K^{-1} \mathcal
J_K^{-\mathsf T} \text{det}(\mathcal J_K)$ in all points. This setup
reduces the memory access besides the source and destination vectors to a
minimum since the merged tensor resides in the cache at all times and thus
better shows the effect of the basis.
We create a mesh
with $2^{l_1}2^{l_2}2^{l_3}$ hexahedral elements and a difference in mesh size of at most
two, i.e., from the mesh sequence $1\times 1\times 1$, $2\times 1\times 1$,
$2\times 2\times 1$, $2\times 2\times 2$, $4\times 2\times 2$, \ldots. The mesh
is selected depending on the polynomial degree to create a problem which
has between $30$ million and $56$ million unknowns. This size makes sure that
the vectors are much larger than the caches of a single node of the two processors
listed in Table~\ref{tab:systems}.  Table~\ref{tab:sizes} lists the meshes chosen for each degree.

\begin{table}
  \caption{Specification of hardware systems used for evaluation with turbo mode enabled on both. Memory
    bandwidth is according to the STREAM triad benchmark (optimized variant without
read for ownership transfer involving two reads and one write) and GFLOP/s are based
on
    the theoretical maximum at the AVX-512 frequency. The \texttt{dgemm} performance is measured for $m=n=k=12{,}000$ with Intel MKL 18.0.2. The processor Xeon
Platinum 8174 is a special model for SuperMUC-NG and not listed in official
Intel documents. We measured a frequency of 2.7 GHz with AVX-512 dense code for
the current experiments. The empirical machine balance is computed as the ratio of measured
\texttt{dgemm} performance and STREAM bandwidth from RAM memory.}
    \label{tab:systems}
  {
	\footnotesize
	\strut\hfill
    \begin{tabular}{lcc}
      \hline
      & Intel Cascade Lake  & Intel Skylake \\
      & Xeon Gold 6230 & Xeon Platinum 8174 \\
      \hline
      cores       & $2\times 20$ & $2\times 24$ \\
      frequency base & 2.1 GHz & 2.7 GHz \\
      max AVX-512 frequency & 2.0 GHz & 2.7 GHz \\
      SIMD width    & 512 bit & 512 bit \\
      arithmetic peak  & 2560 GFLOP/s & 4147 GFLOP/s\\
      \texttt{dgemm} performance & 2124 GFLOP/s & 2920 GFLOP/s \\
      memory interface & DDR4-2933, 12 channels & DDR4-2666, 12 channels \\
      STREAM memory bandwidth & 181 GB/s & 205 GB/s\\
      empirical machine balance & 11.7 FLOP/Byte & 14.3 FLOP/Byte  \\
      \hline
      compiler & \multicolumn{2}{c}{\texttt{g++}, version 9.1.0}\\
      compiler flags & \multicolumn{2}{c}{\texttt{-O3 -funroll-loops -march=skylake-avx512}}\\
      \hline
    \end{tabular}
    \hfill\strut
  }
\end{table}

\begin{table}
  \caption{Number of elements $2^{l_1}2^{l_2}2^{l_3}$ and unknowns for different polynomial degrees.}
  \label{tab:sizes}
  \footnotesize
  \setlength{\tabcolsep}{2pt}
  \centering
  \begin{tabular}{lcccccccccccccccc}
  \hline
  degree $p$ & 1 & 2 & 3 & 4 & 5 & 6 & 7 & 8 & 9 & 10 & 11 & 12 & 13 & 14 & 15 & 16\\
  $l_1+l_2+l_3$ & 22 & 21 & 19 & 18 & 18 & 17 & 16 & 16 & 15 & 15 & 15 & 14 & 14 & 14 & 13 & 13\\
  $n_\text{DoF}$ & 34M & 57M & 34M & 33M & 57M & 45M & 34M & 48M & 33M & 44M & 57M & 36M & 45M & 55M & 34M & 40M\\
  \hline
  \end{tabular}
\end{table}

The number of arithmetic operations per unknown for a matrix-vector product in the given setup is listed in
Table~\ref{tab:operations} for the proposed Hermite-like basis as well as two
nodal bases, a generic one based on the Gauss--Lobatto node points as well as
the collocation case where the nodal points coincide with the quadrature points
and $S=I$ in equations~\eqref{eq:tensor_cell} and \eqref{eq:tensor_face}, which
reduces the number of arithmetic operations. All sum factorization sweeps make
use of the even-odd decomposition \cite{Solomonoff92}.
The ideal data access is
24~Bytes per DoF (3~doubles per DoF) for bases, one double to load the input vector, one
double to write to the output vector, and an additional load operation on the output
vector due to the read-for-ownership (RFO) data transfer \cite{Hager11}. Thus, the
arithmetic intensity is between 8 and 14~FLOP/Byte for the Hermite-like basis
and between 7 and 9.5~FLOP/Byte for the collocated nodal basis. Compared to the
empirical machine balance of 11.7~FLOP/Byte of the Xeon Gold and 14.3~FLOP/Byte
of Xeon Platinum, the roofline model~\cite{williams09} suggests that
the algorithm is in a regime where
both memory transfer and arithmetic matter. The particular instruction mix
with the even-odd decomposition and operations at quadrature points yields
an achievable floating point throughput of around 60--75\% of peak for
$p\leq 10$ due to the proportion of FMA instructions among all arithmetic
instructions \cite[Figure 3]{Kronbichler19}. The effective
machine balance on the Xeon Gold is hence around 7--8~FLOP/Byte and 10--11~FLOP/Byte on the Xeon
Platinum. Assuming ideal execution, we expect the Xeon Gold to be core limited
and the Xeon Platinum to be primarily memory access limited.

\begin{table}
  \caption{Number of arithmetic operations per degree of freedom (FLOP/DoF) for
evaluating the 3D Laplacian with element-wise face integrals for the
Hermite-like basis and nodal bases in Gauss--Lobatto points and Gauss points,
respectively. Gaussian quadrature with $(p+1)^3$ points is used.  The numbers are
based on counting the number of 1D interpolations times the cost of one interpolation per element, plus the work at quadrature points, and dividing by the unknowns per element.}

    \label{tab:operations}
  {
	\footnotesize
	\strut\hfill
    \begin{tabular}{lccccccccc}
      \hline
      polynomial degree $p$ & 1 & 2 & 3 & 4 & 5 & 7 & 9 & 11 & 16\\
      \hline
      Hermite-like basis & 244 & 191 & 218 & 206 & 225 & 241 & 260 & 281 & 333\\
      nodal Gauss--Lobatto basis & 258 & 210 & 240 & 229 & 250 & 267 & 287 & 308 & 361 \\
      nodal Gauss basis & 204 & 168 & 180 & 171 & 180 & 186 & 194 & 204 & 229 \\
      \hline
    \end{tabular}
    \hfill\strut
  }
\end{table}

Our implementation\footnote{The code for all experiments is available on \url{https://github.com/kronbichler/multigrid}, subfolder \texttt{matvec\_dg} for section \ref{sec:experiments} and \texttt{poisson\_dg\_plain}, \texttt{matvec\_dg\_cheby} for section \ref{sec:multigrid}, respectively.} is based on the \texttt{deal.II} finite element library
\cite{dealII90} and traverses the elements in the mesh in the Morton order
(Z-order). In all experiments, vectorization over several elements is chosen
according to \cite{Kronbichler19}, which yields batches of eight elements for double
precision arithmetic and sixteen elements for single-precision arithmetic.
For the given structured meshes, the Morton cell traversal
ensures that a brick of closely packed elements are batched together. Degrees
of freedom are arranged cell-wise for easier implementation of
multigrid solvers and mixed-precision algorithms using the same unknown numbering.
Due to vectorization over cells, this involves transpose
operations within SIMD lanes (array-of-struct to struct-of-array) at the
beginning and the end of the cell access. We use a
pure MPI parallelization with as many ranks as there are cores in the system (no
hyperthreading), with parallel partitioning of the elements created by the
\texttt{p4est} library via the Morton curve \cite{Bangerth11,Burstedde11}.
The evaluation of face integrals requires exchange of data on elements which have
a neighbor that is owned by another process. For the Hermite-like basis, we
only need to exchange the $2(p+1)^{d-1}$ data items per face with ghost neighbor
as given by the access highlighted in Figure~\ref{fig:data_access},
whereas the full ghosted elements are exchanged for
the two nodal bases. In addition, a pure shared-memory parallelization based on OpenMP
is considered, where the loop over element batches is split statically into as
many partitions as there are threads resulting in a similar portioning as in the
MPI case. Threads are pinned to cores with the ``close'' affinity rule.

\pgfplotstableread{
degree ndofs timehermite timegll  timegauss   opsH   opsGL  opsG
1  33554432  0.0133535  0.0141137  0.0121075  244    258    204
2  56623104  0.0175003  0.0180336  0.0182774  190.7  210    168
3  33554432  0.010198   0.010731   0.00930343 218    240    180
4  32768000  0.00951992 0.0106821  0.00917988 205.6  229.2  170.9
5  56623104  0.0156507  0.0178947  0.015893   225.3  250    180
6  44957696  0.0119584  0.0142153  0.0126308  223.8  249.2  177.8
7  33554432  0.00973989 0.0113545  0.00976275 241    267    186
8  47775744  0.0129081  0.015472   0.0137215  244    270.4  186.7
9  32768000  0.00896466 0.010988   0.00962389 260    286.8  194.4
10 43614208  0.0122683  0.0152719  0.0134266  265.4  292.5  196.6
11 56623104  0.0168931  0.0210784  0.0182841  280.6  308    204
12 35995648  0.0111159  0.0139618  0.0121825  287.6  315.1  207.1
13 44957696  0.0145953  0.0183922  0.0165033  302.3  330    214.3
14 55296000  0.0184036  0.0241685  0.0220475  310.1  338    217.9
15 33554432  0.0136616  0.0182003  0.0160831  324.5  352.5  225
16 40247296  0.0156945  0.0231436  0.0220926  333.0  361.1  229.0
}\tableCSLOMP
\pgfplotstableread{
degree ndofs timehermite timegll  timegauss  opsH   opsGL  opsG
1  33554432  0.0148084 0.0161985  0.0133438  244    258    204
2  56623104  0.0188468 0.0205791  0.0191304  190.7  210    168
3  33554432  0.0109867 0.0129177  0.011669   218    240    180
4  32768000  0.0105209 0.0134871  0.0120857  205.6  229.2  170.9
5  56623104  0.0171599 0.0228293  0.0209245  225.3  250    180
6  44957696  0.0134387 0.0188635  0.0173384  223.8  249.2  177.8
7  33554432  0.0107864 0.0166738  0.0149985  241    267    186
8  47775744  0.0143027 0.0231315  0.0211993  244    270.4  186.7
9  32768000  0.0101124 0.0176046  0.0162322  260    286.8  194.4
10 43614208  0.0134776 0.0240553  0.0221749  265.4  292.5  196.6
11 56623104  0.0186298 0.0328294  0.0298277  280.6  308    204
12 35995648  0.0121073 0.0238597  0.0219623  287.6  315.1  207.1
13 44957696  0.0157196 0.0306754  0.0282571  302.3  330    214.3
14 55296000  0.020077  0.0391785  0.036406   310.1  338    217.9
15 33554432  0.0144026 0.0298375  0.0278526  324.5  352.5  225
16 40247296  0.0165197 0.0360423  0.0342649  333.0  361.1  229.0
}\tableCSLMPI

\begin{figure}
  \begin{tikzpicture}
    \begin{axis}[
      width=0.528\textwidth,
      height=0.35\textwidth,
      title style={font=\footnotesize},
      title={throughput OpenMP},
      xlabel={polynomial degree},
      ylabel={GDoF/s},
      legend columns = 4,
      legend to name=legend:cslthroughput,
      legend cell align={left},
      cycle list name=colorGPL,
      tick label style={font=\scriptsize},
      label style={font=\scriptsize},
      legend style={font=\scriptsize},
      grid,
      semithick,
      ymin=0,ymax=4,
      xmin=0.5,xmax=16.5,
      xtick={2,4,6,8,10,12,14,16},
      ]
          \addplot table [x expr={\thisrowno{0}}, y expr={1e-9*\thisrowno{1}/\thisrowno{2}/1}] {\tableCSLOMP};
          \addlegendentry{Hermite-like};
          \addplot table [x expr={\thisrowno{0}}, y expr={1e-9*\thisrowno{1}/\thisrowno{3}/1}] {\tableCSLOMP};
          \addlegendentry{nodal Gauss--Lobatto};
          \addplot table [x expr={\thisrowno{0}}, y expr={1e-9*\thisrowno{1}/\thisrowno{4}/1}] {\tableCSLOMP};
          \addlegendentry{nodal Gauss};
    \end{axis}
  \end{tikzpicture}
  \hfill
  \begin{tikzpicture}
    \begin{axis}[
      width=0.528\textwidth,
      height=0.35\textwidth,
      title style={font=\footnotesize},
      title={throughput MPI},
      xlabel={polynomial degree},
      ylabel={GDoF/s},
      legend pos={south west},
      legend cell align={left},
      cycle list name=colorGPL,
      tick label style={font=\scriptsize},
      label style={font=\scriptsize},
      legend style={font=\scriptsize},
      grid,
      semithick,
      ymin=0,ymax=4,
      xmin=0.5,xmax=16.5,
      xtick={2,4,6,8,10,12,14,16},
      ]
          \addplot table [x expr={\thisrowno{0}}, y expr={1e-9*\thisrowno{1}/\thisrowno{2}/1}] {\tableCSLMPI};
          \addplot table [x expr={\thisrowno{0}}, y expr={1e-9*\thisrowno{1}/\thisrowno{3}/1}] {\tableCSLMPI};
          \addplot table [x expr={\thisrowno{0}}, y expr={1e-9*\thisrowno{1}/\thisrowno{4}/1}] {\tableCSLMPI};
    \end{axis}
  \end{tikzpicture}
  \\
  \begin{tikzpicture}
    \begin{axis}[
      width=0.49\textwidth,
      height=0.35\textwidth,
      title style={font=\footnotesize},
      xlabel={polynomial degree},
      ylabel={GFLOP/s},
      legend pos={south west},
      legend cell align={left},
      cycle list name=colorGPL,
      tick label style={font=\scriptsize},
      label style={font=\scriptsize},
      legend style={font=\scriptsize},
      grid,
      semithick,
      ymin=0,ymax=1000,
      xmin=0.5,xmax=16.5,
      xtick={2,4,6,8,10,12,14,16},
      ]
          \addplot table [x expr={\thisrowno{0}}, y expr={1e-9*\thisrowno{1}/\thisrowno{2}/1*\thisrowno{5}}] {\tableCSLOMP};
          \addplot table [x expr={\thisrowno{0}}, y expr={1e-9*\thisrowno{1}/\thisrowno{3}/1*\thisrowno{6}}] {\tableCSLOMP};
          \addplot table [x expr={\thisrowno{0}}, y expr={1e-9*\thisrowno{1}/\thisrowno{4}/1*\thisrowno{7}}] {\tableCSLOMP};
    \end{axis}
  \end{tikzpicture}
  \hfill
  \begin{tikzpicture}
    \begin{axis}[
      width=0.49\textwidth,
      height=0.35\textwidth,
      title style={font=\footnotesize},
      xlabel={polynomial degree},
      ylabel={GFLOP/s},
      legend pos={south west},
      legend cell align={left},
      cycle list name=colorGPL,
      tick label style={font=\scriptsize},
      label style={font=\scriptsize},
      legend style={font=\scriptsize},
      grid,
      semithick,
      ymin=0,ymax=1000,
      xmin=0.5,xmax=16.5,
      xtick={2,4,6,8,10,12,14,16},
      ]
          \addplot table [x expr={\thisrowno{0}}, y expr={1e-9*\thisrowno{1}/\thisrowno{2}/1*\thisrowno{5}}] {\tableCSLMPI};
          \addplot table [x expr={\thisrowno{0}}, y expr={1e-9*\thisrowno{1}/\thisrowno{3}/1*\thisrowno{6}}] {\tableCSLMPI};
          \addplot table [x expr={\thisrowno{0}}, y expr={1e-9*\thisrowno{1}/\thisrowno{4}/1*\thisrowno{7}}] {\tableCSLMPI};
    \end{axis}
  \end{tikzpicture}
  \\
  \strut\hfill\pgfplotslegendfromname{legend:cslthroughput}\hfill\strut
\caption{Throughput of double-precision matrix-vector product for the 3D
Laplacian on an affine geometry in billion DoF/s for various bases with OpenMP
(left panels) and MPI parallelization (right panels), respectively,
on $2\times 20$ Xeon Gold cores.}
\label{fig:matvec_csl}
\end{figure}
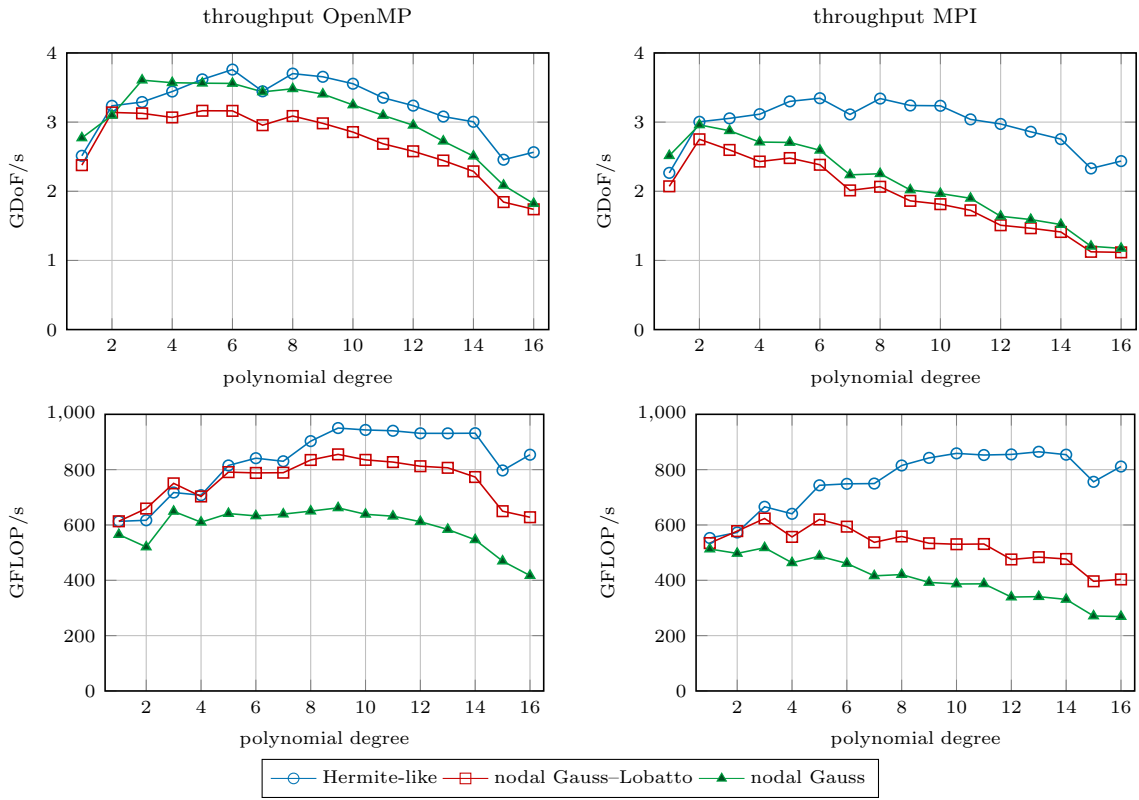

For the performance measurements, we record timings $T$ for $n_\text{tests} =
200$ matrix-vector products for $p=1,\ldots,16$ and report the throughput
\begin{equation}\label{eq:throughput}
  \text{billion DoF / s} = \text{GDoF / s} = 10^{-9}\frac{n_\text{tests} n_\text{DoF}}{T},
\end{equation}
as common in this context~\cite{Fischer19}. We repeat the experiment 40 times
and report the maximum throughput among those experiments. The standard deviation
is 1--3\%, mostly caused by background work of the operating system and the
dynamic frequency adjustment of the processors.
Figure~\ref{fig:matvec_csl} reports the results with the three bases on the
Intel Xeon Gold Cascade Lake for the OpenMP
and MPI parallelizations, respectively. The top two panels display the
application throughput, whereas the bottom two panels the achieved floating
point performance computed from the throughput in DoF/s multiplied by the
work per unknown from Table~\ref{tab:operations}. The OpenMP results
show a similar arithmetic performance for the Hermite-like basis and the nodal
Gauss--Lobatto basis for degrees $p\leq 7$, indicating that both are primarily
core-limited. Due to the cheaper face-normal interpolation via
equation~\eqref{eq:tensor_face}, the Hermite-like basis reaches a higher
DoF/s throughput. The collocated nodal basis using Gauss points with less arithmetic work
can initially outperform the other two bases in terms of GDoF/s. For
higher polynomial degrees, $p\geq 5$, the proposed Hermite-like basis
is the fastest variant with up to 20\% advantage over the nodal Gauss basis.
This result contradicts intuition in spectral element methods where
collocated bases are favored, and underlines the importance of
data access, even on a system like the Xeon Gold that tends to be core-limited.

When turning to the performance with the MPI parallelization, the advantage of
the proposed Hermite-like basis becomes much more significant, achieving e.g.~48\%
higher throughput than the nodal Gauss basis with $p=8$. While the new
basis is able to sustain an approximately constant DoF/s throughput analogously
to the OpenMP parallelization, the throughput with the nodal bases clearly
drops due to the cost of the MPI data exchange as explained below.

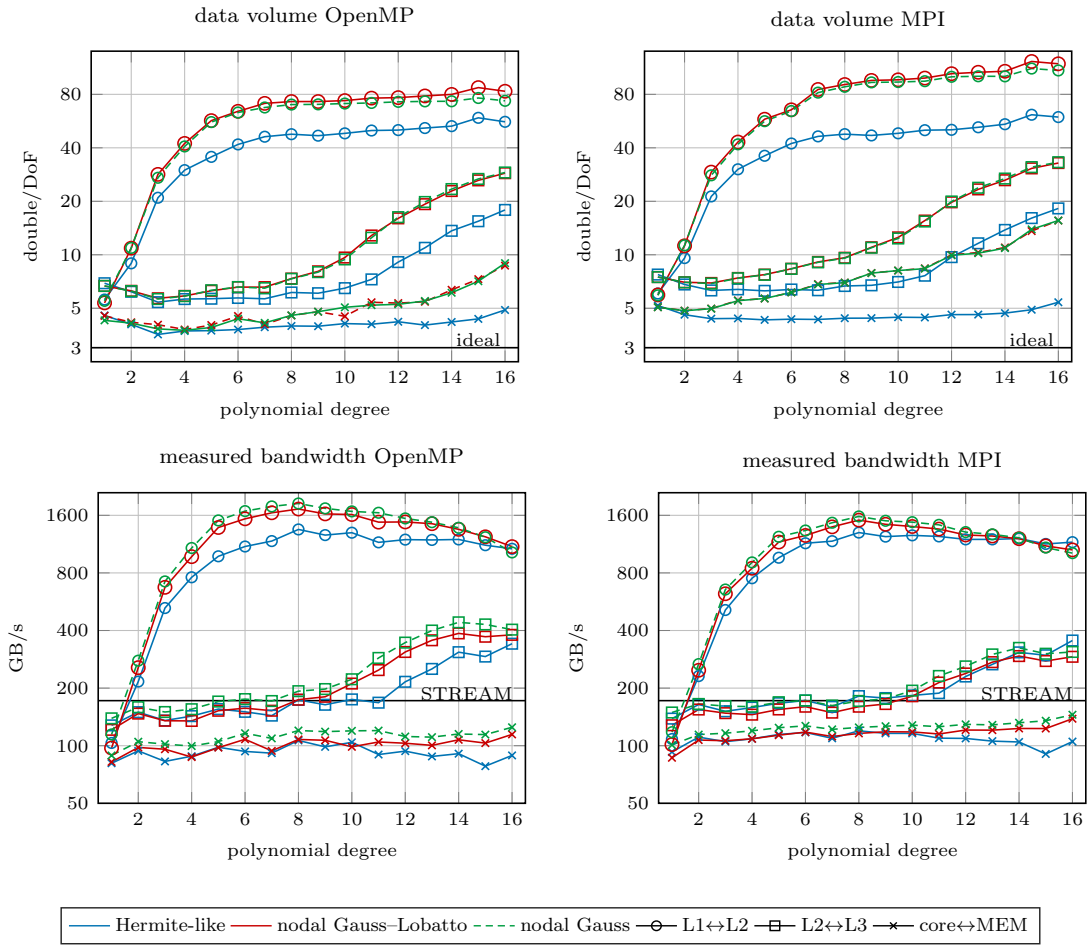
\begin{figure}

  \strut\hfill
  \begin{tikzpicture}
    \begin{semilogyaxis}[
      width=0.48\textwidth,
      height=0.38\textwidth,
      title style={font=\footnotesize},
      title={data volume OpenMP},
      xlabel={polynomial degree},
      ylabel={double/DoF},
      legend pos={south west},
      legend cell align={left},
      cycle list name=colorGPL,
      tick label style={font=\scriptsize},
      label style={font=\scriptsize},
      legend style={font=\scriptsize},
      grid,
      semithick,
      ymin=2.5,ymax=140,
      ytick={3,5,10,20,40,80,160},
      yticklabels={3,5,10,20,40,80,160},
      xmin=0.5,xmax=16.5,
      xtick={2,4,6,8,10,12,14,16},
      ]

	\node[anchor=east] at (axis cs: 0, 3) (nodeAE) {};
	\node[anchor=west] at (axis cs: 17, 3) (nodeAF) {};
	\draw [-] (nodeAE) -- (nodeAF);
      \node[anchor=south] at (axis cs:  15, 2.7) (nodeBE) {{\scriptsize ideal}};

      \addplot +[gnuplot@darkblue,every mark/.append style={fill=gnuplot@darkblue!80!black},mark=o] table [x=deg, y=da1] {\tableCacheOMPH};
      \addplot +[gnuplot@darkblue,every mark/.append style={fill=gnuplot@darkblue!80!black},mark=square] table [x=deg, y=da2] {\tableCacheOMPH};
      \addplot +[gnuplot@darkblue,every mark/.append style={fill=gnuplot@darkblue!80!black},mark=x] table [x=deg, y=da3] {\tableCacheOMPH};

      \addplot +[gnuplot@red,every mark/.append style={fill=gnuplot@red!50!black},mark=o] table [x=deg, y=da1] {\tableCacheOMPL};
      \addplot +[gnuplot@red,every mark/.append style={fill=gnuplot@red!50!black},mark=square] table [x=deg, y=da2] {\tableCacheOMPL};
      \addplot +[gnuplot@red,every mark/.append style={fill=gnuplot@red!50!black},mark=x] table [x=deg, y=da3] {\tableCacheOMPL};

      \addplot +[gnuplot@green,every mark/.append style={fill=gnuplot@green!50!black},mark=o] table [x=deg, y=da1] {\tableCacheOMPG};
      \addplot +[gnuplot@green,every mark/.append style={fill=gnuplot@green!50!black},mark=square] table [x=deg, y=da2] {\tableCacheOMPG};
      \addplot +[gnuplot@green,every mark/.append style={fill=gnuplot@green!50!black},mark=x] table [x=deg, y=da3] {\tableCacheOMPG};

    \end{semilogyaxis}
  \end{tikzpicture}
  \hfill
  \begin{tikzpicture}
    \begin{semilogyaxis}[
      width=0.48\textwidth,
      height=0.38\textwidth,
      title style={font=\footnotesize},
      title={data volume MPI},
      xlabel={polynomial degree},
      ylabel={double/DoF},
      legend pos={south west},
      legend cell align={left},
      cycle list name=colorGPL,
      tick label style={font=\scriptsize},
      label style={font=\scriptsize},
      legend style={font=\scriptsize},
      grid,
      semithick,
      ytick={3,5,10,20,40,80,160},
      ymin=2.5,ymax=140,
      yticklabels={3,5,10,20,40,80,160},
      xmin=0.5,xmax=16.5,
      xtick={2,4,6,8,10,12,14,16},
      ]

	\node[anchor=east] at (axis cs: 0, 3) (nodeAE) {};
	\node[anchor=west] at (axis cs: 17, 3) (nodeAF) {};
	\draw [-] (nodeAE) -- (nodeAF);
    \node[anchor=south] at (axis cs:  15, 2.7) (nodeBE) {{\scriptsize ideal}};

      \addplot +[gnuplot@darkblue,every mark/.append style={fill=gnuplot@darkblue!80!black},mark=o] table [x=deg, y=da1] {\tableCacheMPIH};
      \addplot +[gnuplot@darkblue,every mark/.append style={fill=gnuplot@darkblue!80!black},mark=square] table [x=deg, y=da2] {\tableCacheMPIH};
      \addplot +[gnuplot@darkblue,every mark/.append style={fill=gnuplot@darkblue!80!black},mark=x] table [x=deg, y=da3] {\tableCacheMPIH};

      \addplot +[gnuplot@red,every mark/.append style={fill=gnuplot@red!50!black},mark=o] table [x=deg, y=da1] {\tableCacheMPIL};
      \addplot +[gnuplot@red,every mark/.append style={fill=gnuplot@red!50!black},mark=square] table [x=deg, y=da2] {\tableCacheMPIL};
      \addplot +[gnuplot@red,every mark/.append style={fill=gnuplot@red!50!black},mark=x] table [x=deg, y=da3] {\tableCacheMPIL};

      \addplot +[gnuplot@green,every mark/.append style={fill=gnuplot@green!50!black},mark=o] table [x=deg, y=da1] {\tableCacheMPIG};
      \addplot +[gnuplot@green,every mark/.append style={fill=gnuplot@green!50!black},mark=square] table [x=deg, y=da2] {\tableCacheMPIG};
      \addplot +[gnuplot@green,every mark/.append style={fill=gnuplot@green!50!black},mark=x] table [x=deg, y=da3] {\tableCacheMPIG};

    \end{semilogyaxis}
  \end{tikzpicture}  \hfill\strut

  \strut\hfill
  \begin{tikzpicture}
    \begin{semilogyaxis}[
      width=0.48\textwidth,
      height=0.38\textwidth,
      title style={font=\footnotesize},
      title={measured bandwidth OpenMP},
      xlabel={polynomial degree},
      ylabel={GB/s},
      legend pos={south west},
      legend cell align={left},
      cycle list name=colorGPL,
      tick label style={font=\scriptsize},
      label style={font=\scriptsize},
      legend style={font=\scriptsize},
      grid,
      semithick,
      ymin=50,ymax=2100,
      ytick={50,100,200,400,800,1600},
      yticklabels={50,100,200,400,800,1600},
      xmin=0.5,xmax=16.5,
      xtick={2,4,6,8,10,12,14,16},
      ]

	\node[anchor=east] at (axis cs: 0, 172) (nodeAE) {};
	\node[anchor=west] at (axis cs: 17, 172) (nodeAF) {};
	\draw [-] (nodeAE) -- (nodeAF);
    \node[anchor=south] at (axis cs:  14.3, 155) (nodeBE) {{\scriptsize STREAM}};

      \addplot +[gnuplot@darkblue,solid,every mark/.append style={fill=gnuplot@darkblue!80!black},mark=o] table [x=deg, y=bw1] {\tableCacheOMPH};
      \addplot +[gnuplot@darkblue,solid,every mark/.append style={fill=gnuplot@darkblue!80!black},mark=square] table [x=deg, y=bw2] {\tableCacheOMPH};
      \addplot +[gnuplot@darkblue,solid,every mark/.append style={fill=gnuplot@darkblue!80!black},mark=x] table [x=deg, y=bw3] {\tableCacheOMPH};

      \addplot +[gnuplot@red,solid,every mark/.append style={fill=gnuplot@red!50!black},mark=o] table [x=deg, y=bw1] {\tableCacheOMPL};
      \addplot +[gnuplot@red,solid,every mark/.append style={fill=gnuplot@red!50!black},mark=square] table [x=deg, y=bw2] {\tableCacheOMPL};
      \addplot +[gnuplot@red,solid,every mark/.append style={fill=gnuplot@red!50!black},mark=x] table [x=deg, y=bw3] {\tableCacheOMPL};

      \addplot +[gnuplot@green,densely dashed,every mark/.append style={fill=gnuplot@green!50!black,solid},mark=o] table [x=deg, y=bw1] {\tableCacheOMPG};
      \addplot +[gnuplot@green,densely dashed,every mark/.append style={fill=gnuplot@green!50!black,solid},mark=square] table [x=deg, y=bw2] {\tableCacheOMPG};
      \addplot +[gnuplot@green,densely dashed,every mark/.append style={fill=gnuplot@green!50!black,solid},mark=x] table [x=deg, y=bw3] {\tableCacheOMPG};

    \end{semilogyaxis}
  \end{tikzpicture}
  \hfill
  \begin{tikzpicture}
    \begin{semilogyaxis}[
      width=0.48\textwidth,
      height=0.38\textwidth,
      title style={font=\footnotesize},
      title={measured bandwidth MPI},
      xlabel={polynomial degree},
      ylabel={GB/s},
      legend pos={south west},
      legend cell align={left},
      cycle list name=colorGPL,
      tick label style={font=\scriptsize},
      label style={font=\scriptsize},
      legend style={font=\scriptsize},
      grid,
      semithick,
      ymin=50,ymax=2100,
      ytick={50,100,200,400,800,1600},
      yticklabels={50,100,200,400,800,1600},
      xmin=0.5,xmax=16.5,
      xtick={2,4,6,8,10,12,14,16},
      ]

	\node[anchor=east] at (axis cs: 0, 172) (nodeAE) {};
	\node[anchor=west] at (axis cs: 17, 172) (nodeAF) {};
	\draw [-] (nodeAE) -- (nodeAF);
    \node[anchor=south] at (axis cs:  14.3, 155) (nodeBE) {{\scriptsize STREAM}};

      \addplot +[gnuplot@darkblue,solid,every mark/.append style={fill=gnuplot@darkblue!80!black},mark=o] table [x=deg, y=bw1] {\tableCacheMPIH};
      \addplot +[gnuplot@darkblue,solid,every mark/.append style={fill=gnuplot@darkblue!80!black},mark=square] table [x=deg, y=bw2] {\tableCacheMPIH};
      \addplot +[gnuplot@darkblue,solid,every mark/.append style={fill=gnuplot@darkblue!80!black},mark=x] table [x=deg, y=bw3] {\tableCacheMPIH};

      \addplot +[gnuplot@red,solid,every mark/.append style={fill=gnuplot@red!50!black},mark=o] table [x=deg, y=bw1] {\tableCacheMPIL};
      \addplot +[gnuplot@red,solid,every mark/.append style={fill=gnuplot@red!50!black},mark=square] table [x=deg, y=bw2] {\tableCacheMPIL};
      \addplot +[gnuplot@red,solid,every mark/.append style={fill=gnuplot@red!50!black},mark=x] table [x=deg, y=bw3] {\tableCacheMPIL};

      \addplot +[gnuplot@green,densely dashed,every mark/.append style={fill=gnuplot@green!50!black,solid},mark=o] table [x=deg, y=bw1] {\tableCacheMPIG};
      \addplot +[gnuplot@green,densely dashed,every mark/.append style={fill=gnuplot@green!50!black,solid},mark=square] table [x=deg, y=bw2] {\tableCacheMPIG};
      \addplot +[gnuplot@green,densely dashed,every mark/.append style={fill=gnuplot@green!50!black,solid},mark=x] table [x=deg, y=bw3] {\tableCacheMPIG};

    \end{semilogyaxis}
  \end{tikzpicture}  \hfill\strut

\begin{center}
\begin{tikzpicture}[scale=0.9]
    \begin{axis}[%
    hide axis,
    legend style={font=\footnotesize},
    xmin=10,
    xmax=50,
    ymin=0,
    ymax=0.4,
    semithick,
    legend style={draw=white!15!black,legend cell align=left},legend columns=-1
    ]
    \addlegendimage{gnuplot@darkblue,every mark/.append style={fill=gnuplot@darkblue!80!black}}
    \addlegendentry{Hermite-like};
    \addlegendimage{gnuplot@red,every mark/.append style={fill=gnuplot@red!50!black}}
    \addlegendentry{nodal Gauss--Lobatto};
    \addlegendimage{gnuplot@green,densely dashed,every mark/.append style={fill=gnuplot@green!50!black}}
    \addlegendentry{nodal Gauss};
    \addlegendimage{black,mark=o}
    \addlegendentry{L1$\leftrightarrow$L2};
    \addlegendimage{black,mark=square}
    \addlegendentry{L2$\leftrightarrow$L3};
    \addlegendimage{black,mark=x}
    \addlegendentry{core$\leftrightarrow$MEM};
    \end{axis}
\end{tikzpicture}
\end{center}

\caption{Analysis of data transfer over various memory hierarchies for executing
the 3D DG-SIP Laplacian on $2\times 20$ Xeon Gold cores.}
\label{fig:matvec_cache}
\end{figure}

Further insight into the advantages of the Hermite-like basis is given by an
analysis of the data transfer of the complete matrix-vector product
over the various levels of the memory hierarchy on
the Xeon Gold in Figure~\ref{fig:matvec_cache}. The analysis is based on
hardware performance counters measured with the LIKWID tool \cite{likwid}.
The transfer between the core
and the main memory, indicated by lines with cross marks in the figure, should
ideally amount to 3 doubles per DoF. The lowest recorded values are around 3.5
to 4 doubles per DoF. The excess transfer appears because not all neighboring
data can be cached for elements at far temporal distance (given by the Morton
traversal), an effect well-studied in the context of finite difference stencils
\cite{Hager11}. The measured data transfer from main memory is around
100--120 GB/s and
thus below the STREAM bandwidth, indicating that the execution within the core
puts the primary limit on execution.
However, we see that the nodal bases start to deviate from the
ideal behavior of low degrees once $p>5$. This indicates that access to
all DoF of an element, rather than only two layers with the Hermite-like basis,
is of significance and reduces the effectiveness of caches. The increase of the
memory transfer in the case of the
Hermite-like basis for $p=15,16$ indicates the threshold where the temporary
scratch arrays to perform interpolations of sum factorization exceed the caches,
and the cross-element vectorization should be replaced by intra-element
vectorization \cite{Kronbichler19}. The more beneficial data access of the
proposed Hermite-like basis is also visible for the transfer between the L1 and
L2 cache as well as the transfer between the L2 and L3 cache. Even though the
difference is not tremendous for the OpenMP parallelization, the proposed basis
allows to fit data of around two degrees higher in the same cache level compared
to the nodal bases.

\pgfplotstableread{
deg bw1 bw2 bw3 da1 da2 da3 flops
1.000000000000000000e+00 1.033501061999999990e+02 1.276922887000000060e+02 8.095357889999999657e+01 5.597358942031860352e+00 6.916111335158348083e+00 4.555791243910789490e+00 6.350470344999999952e+02
2.000000000000000000e+00 2.163814442999999983e+02 1.497647318000000212e+02 9.396317690000000766e+01 8.949700108280890021e+00 6.193615330590141355e+00 4.058152878725969970e+00 6.916085357999999133e+02
3.000000000000000000e+00 5.241882891999999856e+02 1.354335217999999941e+02 8.274854039999999600e+01 2.094432674348354340e+01 5.411260947585105896e+00 3.562045469880104065e+00 6.941549473000000035e+02
4.000000000000000000e+00 7.591226944000000003e+02 1.427708245999999974e+02 8.813906080000000998e+01 2.995859794616699290e+01 5.634022140502929155e+00 3.724692916870116832e+00 6.749165454999999838e+02
5.000000000000000000e+00 9.770764392999999473e+02 1.554206345999999996e+02 9.815156500000000506e+01 3.557446334097120655e+01 5.658605142875953398e+00 3.742897289770620972e+00 7.739295413000000963e+02
6.000000000000000000e+00 1.099277799499999901e+03 1.503573200000000156e+02 9.334034440000000643e+01 4.180335915123407631e+01 5.716330803518045123e+00 3.795195598991550234e+00 7.708001045999999405e+02
7.000000000000000000e+00 1.173501569900000050e+03 1.436496625999999992e+02 9.143196490000001120e+01 4.619565233588218689e+01 5.654320865869522095e+00 3.912603482604026794e+00 7.974739205999999285e+02
8.000000000000000000e+00 1.349897595599999931e+03 1.735596587999999940e+02 1.064819595999999962e+02 4.772508016829628730e+01 6.135550971220877337e+00 3.974791318791393469e+00 8.544175699999999551e+02
9.000000000000000000e+00 1.261029278100000056e+03 1.637168918000000133e+02 9.862926530000000014e+01 4.681426086425781818e+01 6.077161788940429688e+00 3.962131881713867099e+00 8.983666853999999375e+02
1.000000000000000000e+01 1.297947825999999850e+03 1.745749898999999914e+02 1.041695869999999928e+02 4.824649658890974280e+01 6.488746453907864087e+00 4.099229544647468693e+00 8.835139307000000599e+02
1.100000000000000000e+01 1.155507885600000009e+03 1.678448086999999873e+02 8.992416770000001236e+01 5.005569480083607203e+01 7.268964582019382448e+00 4.065938569881297759e+00 8.027441593000000921e+02
1.200000000000000000e+01 1.194203988299999992e+03 2.158156448000000012e+02 9.366898609999999792e+01 5.033655283827644666e+01 9.096628153492332558e+00 4.196986535705650390e+00 8.562643901999999798e+02
1.300000000000000000e+01 1.189648883300000080e+03 2.519383533999999827e+02 8.787985760000000823e+01 5.168681825910295657e+01 1.094612616491734869e+01 4.019651329552119456e+00 8.701705093000000488e+02
1.400000000000000000e+01 1.196636565799999971e+03 3.082210090000000378e+02 9.085655130000000668e+01 5.296883997034143476e+01 1.364491215458622797e+01 4.185191062644675775e+00 8.655697787999999946e+02
1.500000000000000000e+01 1.116620852999999897e+03 2.923345876999999859e+02 7.818953700000000140e+01 5.901875793933868408e+01 1.544995978474617004e+01 4.356621578335762024e+00 7.223149331000000757e+02
1.600000000000000000e+01 1.069264100500000040e+03 3.412516865000000053e+02 8.895907090000000039e+01 5.603819912522818214e+01 1.788378404352928186e+01 4.889879123804988126e+00 7.919574516000000131e+02
}\tableCacheOMPH
\pgfplotstableread{
deg bw1 bw2 bw3 da1 da2 da3 flops
1.000000000000000000e+00 9.723344779999999332e+01 1.210519091000000031e+02 8.228123530000000585e+01 5.359264463186264038e+00 6.672076880931854248e+00 4.534758254885673523e+00 6.301128333000000339e+02
2.000000000000000000e+00 2.558990897999999845e+02 1.471842905000000030e+02 9.792500509999999281e+01 1.089410031283343017e+01 6.265908258932607389e+00 4.168565847255566226e+00 7.245818163999999797e+02
3.000000000000000000e+00 6.725079395000000204e+02 1.350189996999999948e+02 9.584865150000000256e+01 2.835034579038619995e+01 5.691881477832794189e+00 4.040131345391273499e+00 7.204117314999999735e+02
4.000000000000000000e+00 9.726855137000000013e+02 1.347635779999999954e+02 8.720946299999999951e+01 4.234873771667480469e+01 5.867332077026366832e+00 3.796498489379882990e+00 6.757299562999999125e+02
5.000000000000000000e+00 1.382063148400000046e+03 1.523073852999999929e+02 9.768821099999999547e+01 5.691502911073190063e+01 6.272201626389114715e+00 4.022620563153867401e+00 7.551768892000000051e+02
6.000000000000000000e+00 1.533818407600000000e+03 1.571621406000000150e+02 1.080372399000000030e+02 6.420198668543868337e+01 6.578432644768984261e+00 4.521814240658596162e+00 7.465134352000000035e+02
7.000000000000000000e+00 1.649157337200000029e+03 1.523719968999999992e+02 9.397754499999999211e+01 7.116097621619701385e+01 6.574836373329162598e+00 4.054661467671394348e+00 7.614868009000000484e+02
8.000000000000000000e+00 1.725837380799999892e+03 1.740605536999999856e+02 1.078359806999999932e+02 7.287316734408155128e+01 7.349674920813373369e+00 4.552998954867138792e+00 7.925180128000000650e+02
9.000000000000000000e+00 1.627329691400000002e+03 1.802213034000000107e+02 1.067398007000000035e+02 7.287317543029784872e+01 8.070460128784178622e+00 4.779455947875977273e+00 8.125575334000000112e+02
1.000000000000000000e+01 1.615340339600000107e+03 2.112628595999999845e+02 9.866756440000000339e+01 7.378036390113973653e+01 9.649393381166062156e+00 4.506322515360132286e+00 7.900049307000000454e+02
1.100000000000000000e+01 1.473901061300000038e+03 2.484091920000000187e+02 1.045713587999999987e+02 7.618146671189202834e+01 1.283951843226397393e+01 5.404666397306654879e+00 7.677877651000000014e+02
1.200000000000000000e+01 1.477199862899999971e+03 3.089443571000000475e+02 1.030244457000000011e+02 7.667259039203850080e+01 1.603544961879835995e+01 5.346965833202947493e+00 7.477983137000001079e+02
1.300000000000000000e+01 1.446220977600000197e+03 3.556791732000000366e+02 1.004211423999999937e+02 7.836550470691381065e+01 1.927297497629771783e+01 5.441097059333289465e+00 7.378808431000001065e+02
1.400000000000000000e+01 1.351548817000000099e+03 3.863930442999999855e+02 1.073744571000000008e+02 7.994226548936630650e+01 2.285462284794559906e+01 6.350750506365741010e+00 7.126244945999999345e+02
1.500000000000000000e+01 1.235269021500000008e+03 3.712170734000000039e+02 1.031861205000000012e+02 8.730499409139156342e+01 2.623647823929786682e+01 7.292461022734642029e+00 6.042692240000000083e+02
1.600000000000000000e+01 1.097770071199999848e+03 3.796665034999999762e+02 1.145541769999999957e+02 8.326399026906055667e+01 2.879706005591033602e+01 8.688282127077556538e+00 5.859100945000000138e+02
}\tableCacheOMPL
\pgfplotstableread{
deg bw1 bw2 bw3 da1 da2 da3 flops
1.000000000000000000e+00 1.144789609000000041e+02 1.389890352999999834e+02 8.891498829999999032e+01 5.503031238913536072e+00 6.681235134601593018e+00 4.273685440421104431e+00 5.718290693000000147e+02
2.000000000000000000e+00 2.782316207000000077e+02 1.589484096000000193e+02 1.049317682000000076e+02 1.093108057975769043e+01 6.244717262409350766e+00 4.122211535771687529e+00 5.935470232000000124e+02
3.000000000000000000e+00 7.221710967000000210e+02 1.499992733999999928e+02 1.019009673999999990e+02 2.713313437998294830e+01 5.635715276002883911e+00 3.828104957938193831e+00 6.031492352000000210e+02
4.000000000000000000e+00 1.079561090700000022e+03 1.551410958999999821e+02 9.960812609999999268e+01 4.071679763793945028e+01 5.851312637329102451e+00 3.756329727172851740e+00 5.798247297000000344e+02
5.000000000000000000e+00 1.504390117299999929e+03 1.701609066999999982e+02 1.050509974999999940e+02 5.583804338066666162e+01 6.315816331792761140e+00 3.898850193730106994e+00 6.149074865000000045e+02
6.000000000000000000e+00 1.683233552699999791e+03 1.752299682999999959e+02 1.162251070000000084e+02 6.323219082445861261e+01 6.582671629791704504e+00 4.365677180165104865e+00 5.870119435999999951e+02
7.000000000000000000e+00 1.780666521500000044e+03 1.712160248999999794e+02 1.092046245000000084e+02 6.747895292937755585e+01 6.488287076354026794e+00 4.137832671403884888e+00 6.087364902999999003e+02
8.000000000000000000e+00 1.844893179099999998e+03 1.927811081000000115e+02 1.200765841999999992e+02 7.006870238378705551e+01 7.321791001308112001e+00 4.560082047073929168e+00 6.170765890999999783e+02
9.000000000000000000e+00 1.739823817700000063e+03 1.977526818000000048e+02 1.184659328999999985e+02 7.025916938781739418e+01 7.985832977294921697e+00 4.783452224731444602e+00 6.214609675000000379e+02
1.000000000000000000e+01 1.677409797699999899e+03 2.224640824000000237e+02 1.195739413999999954e+02 7.101534034505452553e+01 9.418310714251649785e+00 5.061926436907898896e+00 6.086166087999999945e+02
1.100000000000000000e+01 1.650186563799999931e+03 2.875485887000000389e+02 1.199828246000000007e+02 7.157716949780781590e+01 1.247248163929692133e+01 5.203989258518924466e+00 5.729960257999999840e+02
1.200000000000000000e+01 1.541301089499999989e+03 3.458331023000000073e+02 1.117366818999999936e+02 7.259195056024550752e+01 1.628799160665200318e+01 5.262058693873214565e+00 5.534778970999999501e+02
1.300000000000000000e+01 1.472955427200000031e+03 3.998594911999999795e+02 1.108597158999999976e+02 7.303236163392358549e+01 1.982591173044098909e+01 5.496312577939937327e+00 5.242862006000000292e+02
1.400000000000000000e+01 1.377180530299999873e+03 4.413912930000000188e+02 1.151416969999999935e+02 7.312818965205440236e+01 2.343784903067129477e+01 6.113704201027198692e+00 4.873747433999999998e+02
1.500000000000000000e+01 1.228973208099999965e+03 4.304040810999999849e+02 1.143050230000000056e+02 7.631211280822753906e+01 2.672560438513755798e+01 7.097112014889717102e+00 4.203408172999999692e+02
1.600000000000000000e+01 1.023549216399999978e+03 4.042083812000000194e+02 1.248386938999999956e+02 7.356873148695504483e+01 2.905292662394015224e+01 8.972468610065131855e+00 3.840780713999999989e+02
}\tableCacheOMPG
\pgfplotstableread{
deg bw1 bw2 bw3 da1 da2 da3 flops
1.000000000000000000e+00 1.083447465999999935e+02 1.385462320999999974e+02 9.275307879999999727e+01 6.052008271217346191e+00 7.739027962088585788e+00 5.180737748742103577e+00 6.211751923000000488e+02
2.000000000000000000e+00 2.310662996000000078e+02 1.641271357999999907e+02 1.109333583000000090e+02 9.581830766465927240e+00 6.806005813457348452e+00 4.599966605504353545e+00 6.463091513000000532e+02
3.000000000000000000e+00 5.118195670000000064e+02 1.513022643000000187e+02 1.049472734999999943e+02 2.130774743854999542e+01 6.298920884728431702e+00 4.368615895509719849e+00 6.605811115999999856e+02
4.000000000000000000e+00 7.504136581000000206e+02 1.582206998000000056e+02 1.086350424000000032e+02 3.028374671936035512e+01 6.385165023803710760e+00 4.383651733398436612e+00 6.547025254999999788e+02
5.000000000000000000e+00 9.604324043999999958e+02 1.667963005000000010e+02 1.144959401999999926e+02 3.608298058863039159e+01 6.266455076358936616e+00 4.301300313737657532e+00 7.493945039000000179e+02
6.000000000000000000e+00 1.143206282899999906e+03 1.722745151999999962e+02 1.172842693000000054e+02 4.231385183751409329e+01 6.376448873180690313e+00 4.340731784831678652e+00 7.580029549999999290e+02
7.000000000000000000e+00 1.172221391600000061e+03 1.595858251999999879e+02 1.093031509999999997e+02 4.639059640467166901e+01 6.315599009394645691e+00 4.325156286358833313e+00 7.538241709000000128e+02
8.000000000000000000e+00 1.300034990599999901e+03 1.816699260000000038e+02 1.196676184000000092e+02 4.775374167904114131e+01 6.673217114107108827e+00 4.395498267907664847e+00 8.250484558999999081e+02
9.000000000000000000e+00 1.238050209800000175e+03 1.774434999000000062e+02 1.159971287000000046e+02 4.696147842407226847e+01 6.730751800537110086e+00 4.399469375610351562e+00 8.496380864000000201e+02
1.000000000000000000e+01 1.257092182999999977e+03 1.831403818999999942e+02 1.159579453999999998e+02 4.826972382486000868e+01 7.032208689883809427e+00 4.452110537006656266e+00 8.466122378000000026e+02
1.100000000000000000e+01 1.240938743800000111e+03 1.882468466999999919e+02 1.095222799000000009e+02 5.027586906044571435e+01 7.626705920254742921e+00 4.437031790062233583e+00 8.587386711999999989e+02
1.200000000000000000e+01 1.197441811800000096e+03 2.299725620000000106e+02 1.091591391999999985e+02 5.051367668113656606e+01 9.701315489583629770e+00 4.604523302372553140e+00 8.399650640999999496e+02
1.300000000000000000e+01 1.197285880599999928e+03 2.661403482000000054e+02 1.056756100000000060e+02 5.220026594556802735e+01 1.160340767907679194e+01 4.607095245717218290e+00 8.564170934999999645e+02
1.400000000000000000e+01 1.209994279600000027e+03 3.074456541999999786e+02 1.045885260999999957e+02 5.430849790219907192e+01 1.379916766131365691e+01 4.694046359592014106e+00 8.519662327999999434e+02
1.500000000000000000e+01 1.133421931600000107e+03 2.969056930999999508e+02 9.054750269999999546e+01 6.141509115695954080e+01 1.608799360692501423e+01 4.906197264790534973e+00 7.418232232999999951e+02
1.600000000000000000e+01 1.158420992700000170e+03 3.546100678000000244e+02 1.051111570999999998e+02 5.958306584124309069e+01 1.823926693112501241e+01 5.405753531864599282e+00 7.964831445000000940e+02
}\tableCacheMPIH
\pgfplotstableread{
deg bw1 bw2 bw3 da1 da2 da3 flops
1.000000000000000000e+00 1.013800772000000023e+02 1.278647182999999927e+02 8.649004829999999799e+01 5.941778793931007385e+00 7.494017109274864197e+00 5.068599432706832886e+00 6.312973959999999352e+02
2.000000000000000000e+00 2.478796605000000000e+02 1.545622693999999910e+02 1.070474865999999992e+02 1.123064712241843921e+01 7.002728515201145143e+00 4.849710067113240264e+00 6.521855936000000611e+02
3.000000000000000000e+00 6.226013331000000335e+02 1.479651314999999840e+02 1.061733660999999955e+02 2.919679209589958191e+01 6.938801705837249756e+00 4.978596419095993042e+00 6.457606036999999333e+02
4.000000000000000000e+00 8.461684522000000470e+02 1.453113868000000082e+02 1.086640246999999988e+02 4.312074737548827841e+01 7.405060577392577947e+00 5.537229537963867188e+00 5.752903871999999410e+02
5.000000000000000000e+00 1.159548624499999960e+03 1.547089454000000046e+02 1.134624155000000059e+02 5.796631331796999831e+01 7.733959842611242586e+00 5.671827881424515105e+00 6.238045285999999123e+02
6.000000000000000000e+00 1.255486342700000023e+03 1.599290311000000031e+02 1.178119934000000057e+02 6.547394855376930423e+01 8.340325925510061467e+00 6.143913046611641171e+00 5.982526096999999936e+02
7.000000000000000000e+00 1.392508689100000083e+03 1.488980696999999793e+02 1.116994255999999979e+02 8.522776104509830475e+01 9.113216400146484375e+00 6.836464628577232361e+00 5.378878048999999919e+02
8.000000000000000000e+00 1.509779736500000126e+03 1.599153495999999848e+02 1.159451747000000097e+02 9.085977865671753761e+01 9.623826580283083132e+00 6.977987354838471212e+00 5.548377369999999473e+02
9.000000000000000000e+00 1.435932503300000008e+03 1.651472573000000068e+02 1.182654761000000008e+02 9.582358322143552698e+01 1.102069778442382564e+01 7.892236709594726562e+00 5.301217447000000220e+02
1.000000000000000000e+01 1.398407644699999992e+03 1.813523305000000221e+02 1.180482676000000026e+02 9.659504736392324276e+01 1.252690980884027461e+01 8.154535329404582811e+00 5.202645145000000184e+02
1.100000000000000000e+01 1.356084622699999954e+03 2.138305770999999993e+02 1.151566468000000043e+02 9.879552964810972071e+01 1.557830969492594697e+01 8.390353123346963926e+00 5.225439440000000104e+02
1.200000000000000000e+01 1.262941545500000075e+03 2.384879599999999868e+02 1.207464030999999949e+02 1.045434402236625857e+02 1.974149069354161767e+01 9.995295820205820192e+00 4.670475214999999594e+02
1.300000000000000000e+01 1.247225529499999993e+03 2.733692968000000292e+02 1.205917434000000128e+02 1.064984924939213755e+02 2.334256269493881319e+01 1.029724115532966877e+01 4.781064473000000135e+02
1.400000000000000000e+01 1.207367836699999998e+03 2.941341462000000320e+02 1.230362314000000055e+02 1.078991721824363452e+02 2.628598022460937500e+01 1.099558263708044059e+01 4.660880771999999865e+02
1.500000000000000000e+01 1.107499847999999929e+03 2.770698712999999884e+02 1.229319871000000006e+02 1.224874589592218399e+02 3.064343705773353577e+01 1.359703727066516876e+01 3.922042615000000296e+02
1.600000000000000000e+01 1.055309158899999829e+03 2.916034240999999838e+02 1.383349935999999900e+02 1.188714907704607100e+02 3.284664794623717654e+01 1.558210177647710637e+01 3.938652423999999996e+02
}\tableCacheMPIL
\pgfplotstableread{
deg bw1 bw2 bw3 da1 da2 da3 flops
1.000000000000000000e+00 1.176976873999999924e+02 1.491400983999999994e+02 1.008737670999999949e+02 5.912321433424949646e+00 7.491771504282951355e+00 5.067102238535881042e+00 5.861959564999999657e+02
2.000000000000000000e+00 2.670381313000000318e+02 1.652453150000000051e+02 1.143609242999999935e+02 1.130413942866855059e+01 6.995090511110093168e+00 4.840836701569734224e+00 5.396677971999999954e+02
3.000000000000000000e+00 6.550708227000000079e+02 1.610081951999999887e+02 1.162740534999999937e+02 2.801102101802825928e+01 6.884755194187164307e+00 4.971490800380706787e+00 5.332307550999998966e+02
4.000000000000000000e+00 9.078380529999999453e+02 1.600252693000000193e+02 1.196218673999999993e+02 4.185746231079101420e+01 7.378237152099609730e+00 5.514744567871093572e+00 4.708052004999999554e+02
5.000000000000000000e+00 1.235840634000000136e+03 1.687733575999999971e+02 1.241702456999999953e+02 5.649505831577159398e+01 7.715272903442382812e+00 5.676162463647346890e+00 4.911211050000000000e+02
6.000000000000000000e+00 1.332597746700000016e+03 1.723074043000000017e+02 1.270013993000000028e+02 6.443012865917327758e+01 8.330920138790030904e+00 6.140359706155760655e+00 4.603326145000000338e+02
7.000000000000000000e+00 1.463326199599999882e+03 1.621843682999999885e+02 1.217386521000000101e+02 8.180738389492034912e+01 9.066919237375259399e+00 6.805779412388801575e+00 4.098846705999999926e+02
8.000000000000000000e+00 1.578999395000000050e+03 1.713180070000000228e+02 1.246692766000000034e+02 8.822914641580463524e+01 9.572659140588161364e+00 6.966472358860596792e+00 4.150363896000000068e+02
9.000000000000000000e+00 1.497829494399999930e+03 1.756847396999999944e+02 1.264586318999999861e+02 9.340016098022459801e+01 1.095515594482421839e+01 7.885799026489256569e+00 3.827565759999999955e+02
1.000000000000000000e+01 1.473548914100000047e+03 1.940059258000000000e+02 1.280432681000000059e+02 9.381897121690252561e+01 1.235209292348034715e+01 8.152529962254503104e+00 3.783484381999999755e+02
1.100000000000000000e+01 1.426600923500000135e+03 2.315001911999999891e+02 1.256569812999999982e+02 9.486296088607222998e+01 1.539379159609476844e+01 8.356008043995609924e+00 3.769933035999999902e+02
1.200000000000000000e+01 1.307396340099999861e+03 2.597450092000000268e+02 1.293312962999999911e+02 1.004287348153865764e+02 1.995253565097647552e+01 9.934713149211816585e+00 3.334261501000000294e+02
1.300000000000000000e+01 1.272690089300000182e+03 3.000042271000000369e+02 1.284736871000000065e+02 1.011278388799105841e+02 2.383833181753798058e+01 1.020873016935743216e+01 3.327292400000000043e+02
1.400000000000000000e+01 1.220746290199999976e+03 3.245800975999999878e+02 1.317912197999999933e+02 1.009887349446614593e+02 2.685157696759259238e+01 1.090300993742766344e+01 3.243222132000000215e+02
1.500000000000000000e+01 1.085249458900000036e+03 3.025114194000000225e+02 1.350364332000000047e+02 1.117537245154380798e+02 3.115118108689785004e+01 1.390650831162929535e+01 2.668477687000000174e+02
1.600000000000000000e+01 1.015032239300000015e+03 3.090785268000000201e+02 1.451561509000000001e+02 1.090454421459767360e+02 3.320452043536043618e+01 1.559443272909563838e+01 2.618856741999999826e+02
}\tableCacheMPIG

In the MPI case, the data transfer is significantly higher because the
implementation must pack the
data for transfer into a separate buffer and eventually perform a mem-copy
operation from one process to the other, despite actually running in shared
memory. Adding these operations up, and
considering that the ghost data is almost as big as the locally owned part for
$p=16$, the transfer from main memory is more than doubled. While the
Hermite-like basis involves a main memory transfer of around 4--4.5 double per
DoF (increasing the data access by around 0.5 double/DoF over the OpenMP case), the transfer rises
to 10 double/DoF for the two nodal bases for $p=12$, for instance. We note that
this transfer cannot be overlapped with the local computations as it is not
happening over an Infiniband fabric but either in user code for pack/unpack or a
big \texttt{memcpy} call of the MPI library. As a consequence, the execution
stalls at around 130--140 GB/s with a mix of phases with 100 GB/s (computation
phase---compare to the OpenMP parallelization) and others with
180 GB/s (under the assumption that data exchange is performed with full
bandwidth measured with the STREAM benchmark).

\pgfplotstableread{
degree ndofs timehermite timegll  timegauss   opsH   opsGL  opsG
1  33554432  0.0087389 0.0092676  0.00787527 244    258    204
2  56623104  0.012845  0.0134812  0.0128625  190.7  210    168
3  33554432  0.0066677 0.0073904  0.00677957 218    240    180
4  32768000  0.0064200 0.0077299  0.00717094 205.6  229.2  170.9
5  56623104  0.0107926 0.0131163  0.0117841  225.3  250    180
6  44957696  0.0085192 0.0103738  0.00962975 223.8  249.2  177.8
7  33554432  0.0065618 0.007903   0.00714877 241    267    186
8  47775744  0.0092831 0.0114031  0.0105196  244    270.4  186.7
9  32768000  0.0067124 0.0085155  0.00760561 260    286.8  194.4
10 43614208  0.0094964 0.0122245  0.0110163  265.4  292.5  196.6
11 56623104  0.0130022 0.0177083  0.0161904  280.6  308    204
12 35995648  0.0093514 0.0130861  0.0126261  287.6  315.1  207.1
13 44957696  0.0130422 0.0207161  0.0201924  302.3  330    214.3
14 55296000  0.0184955 0.0315809  0.0316683  310.1  338    217.9
15 33554432  0.0138953 0.0246557  0.0248694  324.5  352.5  225
16 40247296  0.0190564 0.0340302  0.0345184  333.0  361.1  229.0
}\tableSKXOMP
\pgfplotstableread{
degree ndofs timehermite timegll  timegauss  opsH   opsGL  opsG
1  33554432  0.0117982  0.0125377  0.0106682  244    258    204
2  56623104  0.0148359  0.0163751  0.0152357  190.7  210    168
3  33554432  0.0088912  0.010533   0.0096135  218    240    180
4  32768000  0.00855881 0.0113134  0.0101973  205.6  229.2  170.9
5  56623104  0.0139139  0.0186706  0.0172089  225.3  250    180
6  44957696  0.0108579  0.0152094  0.0141846  223.8  249.2  177.8
7  33554432  0.00848726 0.0133078  0.0121764  241    267    186
8  47775744  0.0112478  0.0183545  0.0172006  244    270.4  186.7
9  32768000  0.00782293 0.0142494  0.0134173  260    286.8  194.4
10 43614208  0.010577   0.0194269  0.0180849  265.4  292.5  196.6
11 56623104  0.0142216  0.0260703  0.0239724  280.6  308    204
12 35995648  0.00956939 0.0185817  0.017325   287.6  315.1  207.1
13 44957696  0.0121484  0.0239494  0.0220821  302.3  330    214.3
14 55296000  0.0155314  0.0306122  0.0284754  310.1  338    217.9
15 33554432  0.0112653  0.0231333  0.0215816  324.5  352.5  225
16 40247296  0.0131785  0.0276901  0.0267735  333.0  361.1  229.0
}\tableSKXMPI
\pgfplotstableread{
degree ndofs hermiteTwoH gaussTwoH hermiteThH gaussThH  hermiteTwoM gaussTwoM hermiteThM gaussThM
1  33554432  0.0088786  0.0088552  0.010477   0.0099770 0.0118404  0.0111192  0.0121269 0.0111344
2  56623104  0.0123018  0.0124751  0.015361   0.0156492 0.0149848  0.0153478  0.0153909 0.0159258
3  33554432  0.0076100  0.0074884  0.008645   0.0092861 0.00900469 0.00976468 0.0093268 0.0102538
4  32768000  0.0065879  0.0076428  0.00790657 0.0096068 0.00861692 0.0103002  0.0089390 0.0109582
5  56623104  0.011182   0.0128129  0.0128539  0.0160075 0.0140829  0.017425   0.0145149 0.0185082
6  44957696  0.0091239  0.0108574  0.0105764  0.0140701 0.0109779  0.0145408  0.0114205 0.0153668
7  33554432  0.0069598  0.0083004  0.00822133 0.0113462 0.00858913 0.0123338  0.0088982 0.013327
8  47775744  0.0096200  0.0117429  0.0113618  0.0152648 0.0113364  0.0174926  0.0118039 0.0188409
9  32768000  0.0070318  0.0086837  0.0080615  0.0119475 0.00790117 0.0135426  0.0082866 0.0146811
10 43614208  0.0098506  0.0124568  0.0107586  0.0155203 0.0107801  0.018364   0.0110862 0.019772
11 56623104  0.0136791  0.0177286  0.0148906  0.0221485 0.014419   0.0244554  0.0150583 0.0270832
12 35995648  0.0097613  0.0147658  0.0107167  0.0171442 0.00959841 0.0178106  0.0102989 0.020092
13 44957696  0.0136267  0.0227461  0.0148347  0.0258733 0.0122466  0.0227591  0.014035  0.0266593
14 55296000  0.0192605  0.043664   0.0204472  0.0424444 0.0162546  0.0311986  0.017134  0.0329744
15 33554432  0.0144661  0.0329806  0.0152946  0.0348767 0.0117213  0.0220063  0.0150496 0.0285209
16 40247296  0.0196579  0.0423457  0.0210824  0.0432834 0.013398   0.0276627  0.0176196 0.0321953
}\tableSKXBig

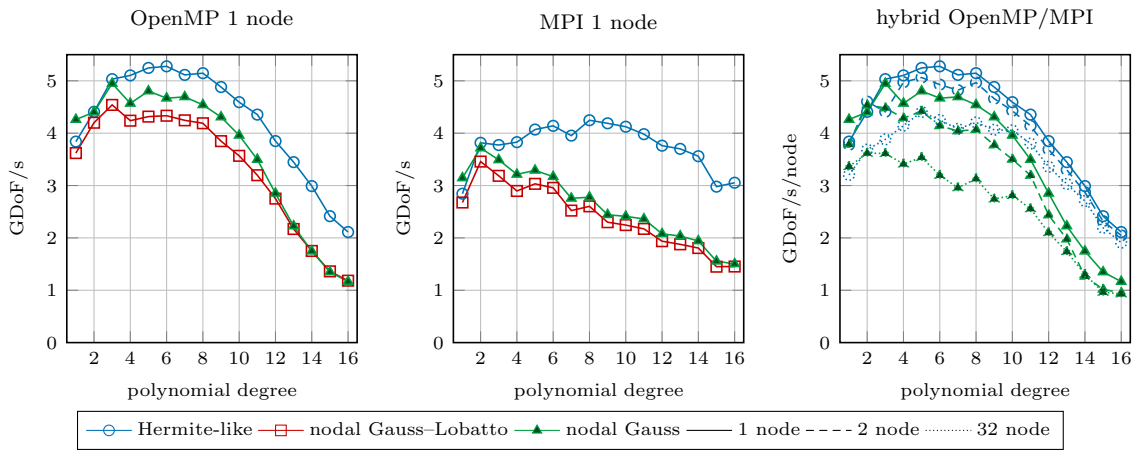
\begin{figure}
  \begin{tikzpicture}
    \begin{axis}[
      width=0.36\textwidth,
      height=0.36\textwidth,
      title style={font=\footnotesize},
      title={OpenMP 1 node},
      xlabel={polynomial degree},
      ylabel={GDoF/s},
      legend columns = 6,
      legend to name=legend:skxthroughput,
      legend cell align={left},
      cycle list name=colorGPL,
      tick label style={font=\scriptsize},
      label style={font=\scriptsize},
      legend style={font=\scriptsize},
      grid,
      semithick,
      ymin=0,ymax=5.5,
      ytick={0,1,2,3,4,5,6},
      xmin=0.5,xmax=16.5,
      xtick={2,4,6,8,10,12,14,16},
      ]
          \addplot table [x expr={\thisrowno{0}}, y expr={1e-9*\thisrowno{1}/\thisrowno{2}/1}] {\tableSKXOMP};
          \addlegendentry{Hermite-like};
          \addplot table [x expr={\thisrowno{0}}, y expr={1e-9*\thisrowno{1}/\thisrowno{3}/1}] {\tableSKXOMP};
          \addlegendentry{nodal Gauss--Lobatto};
          \addplot table [x expr={\thisrowno{0}}, y expr={1e-9*\thisrowno{1}/\thisrowno{4}/1}] {\tableSKXOMP};
          \addlegendentry{nodal Gauss};
          \addplot[black] coordinates {(1,8) (16,8)};
          \addlegendentry{1 node};
          \addplot[black,densely dashed] coordinates {(1,8) (16,8)};
          \addlegendentry{2 node};
          \addplot[black,densely dotted] coordinates {(1,8) (16,8)};
          \addlegendentry{32 node};
    \end{axis}
  \end{tikzpicture}
  \hfill
  \begin{tikzpicture}
    \begin{axis}[
      width=0.36\textwidth,
      height=0.36\textwidth,
      title style={font=\footnotesize},
      title={MPI 1 node},
      xlabel={polynomial degree},
      ylabel={GDoF/s},
      legend pos={south west},
      legend cell align={left},
      cycle list name=colorGPL,
      tick label style={font=\scriptsize},
      label style={font=\scriptsize},
      legend style={font=\scriptsize},
      grid,
      semithick,
      ymin=0,ymax=5.5,
      ytick={0,1,2,3,4,5,6},
      xmin=0.5,xmax=16.5,
      xtick={2,4,6,8,10,12,14,16},
      ]
          \addplot table [x expr={\thisrowno{0}}, y expr={1e-9*\thisrowno{1}/\thisrowno{2}/1}] {\tableSKXMPI};
          \addplot table [x expr={\thisrowno{0}}, y expr={1e-9*\thisrowno{1}/\thisrowno{3}/1}] {\tableSKXMPI};
          \addplot table [x expr={\thisrowno{0}}, y expr={1e-9*\thisrowno{1}/\thisrowno{4}/1}] {\tableSKXMPI};
    \end{axis}
  \end{tikzpicture}
  \hfill
  \begin{tikzpicture}
    \begin{axis}[
      width=0.36\textwidth,
      height=0.36\textwidth,
      title style={font=\footnotesize},
      title={hybrid OpenMP/MPI},
      xlabel={polynomial degree},
      ylabel={GDoF/s/node},
      cycle list name=colorGPL,
      tick label style={font=\scriptsize},
      label style={font=\scriptsize},
      grid,
      semithick,
      ymin=0,ymax=5.5,
      ytick={0,1,2,3,4,5,6},
      xmin=0.5,xmax=16.5,
      xtick={2,4,6,8,10,12,14,16},
      ]
          \addplot[color=gnuplot@darkblue,mark=o] table [x expr={\thisrowno{0}}, y expr={1e-9*\thisrowno{1}/\thisrowno{2}/1}] {\tableSKXOMP};
          \addplot[color=gnuplot@darkblue,densely dashed,mark=o] table [x expr={\thisrowno{0}}, y expr={1e-9*\thisrowno{1}/\thisrowno{2}/1}] {\tableSKXBig};
          \addplot[color=gnuplot@darkblue,densely dotted,mark=o] table [x expr={\thisrowno{0}}, y expr={1e-9*\thisrowno{1}/\thisrowno{4}/1}] {\tableSKXBig};
          \addplot[gnuplot@green,every mark/.append style={fill=gnuplot@green!50!black},mark=triangle*] table [x expr={\thisrowno{0}}, y expr={1e-9*\thisrowno{1}/\thisrowno{4}/1}] {\tableSKXOMP};
          \addplot[gnuplot@green,densely dashed,every mark/.append style={fill=gnuplot@green!50!black},mark=triangle*] table [x expr={\thisrowno{0}}, y expr={1e-9*\thisrowno{1}/\thisrowno{3}/1}] {\tableSKXBig};
          \addplot[gnuplot@green,densely dotted,every mark/.append style={fill=gnuplot@green!50!black},mark=triangle*] table [x expr={\thisrowno{0}}, y expr={1e-9*\thisrowno{1}/\thisrowno{5}/1}] {\tableSKXBig};
    \end{axis}
  \end{tikzpicture}
  \\
  \strut\hfill\pgfplotslegendfromname{legend:skxthroughput}\hfill\strut
\caption{Throughput of double-precision matrix-vector product for the 3D
Laplacian on an affine geometry in billion DoF/s for various bases with OpenMP and hyperthreading
(left panel) and MPI parallelization (middle panel) on 1 node as well as up to 32 nodes with hybrid OpenMP/MPI (right panel), respectively,
on $2\times 24$ Xeon Platinum cores.}
\label{fig:matvec_skx}
\end{figure}

The throughput on the Xeon Platinum is shown in Figure~\ref{fig:matvec_skx}.
Overall, the trend is very similar to the Xeon Gold. As opposed to the
experiment on the Xeon Gold of Figure~\ref{fig:matvec_csl}, the OpenMP
experiment is run with 2-way hyperthreading which gives around 5\% higher
throughput for $1\leq p \leq 10$ than running without. For very
high degrees $p\geq 12$, the caches are strained more with hyperthreading and
twice the number of local scratch arrays for
sum factorization which eventually spill to the slow main memory, such
that throughput is even lower than on the Xeon Gold for $p=15,16$. On the other
hand, the MPI parallelization on the Xeon Platinum does not profit from
hyperthreading even for low degrees, so it is not used. A notable difference
to the Xeon Gold is
that the nodal Gauss--Lobatto and nodal Gauss bases show a more similar
performance on the Xeon Platinum, which is expected from the preliminary
characterization in terms of roofline: While Xeon Gold is core-limited,
the Xeon Platinum is more strongly limited by the available memory bandwidth. At
the same time, the nodal Gauss--Lobatto and nodal Gauss bases have the same
memory access pattern and only differ in the arithmetic work. With respect to
arithmetic performance, we record up to 1.3 TFLOP/s on the Xeon Platinum (for
polynomial degrees $p=8,9$), compared to 0.92 TFLOP/s on the Xeon Gold according
to Figure~\ref{fig:matvec_csl}. The measured memory throughput on the Xeon
Platinum with the Hermite-like basis is up to 185 GB/s, or 127 GB/s in terms of
theoretical throughput with 3 doubles per unknown. For the MPI case, the Xeon
Platinum has less advantage over the Xeon Gold, which can be explained by the
higher proportion of time spent in the memory-limited pack/unpack and MPI transfer routines.

The right panel of Figure~\ref{fig:matvec_skx} compares the throughput on one node with the
throughput on 2 and 32 nodes of Xeon Platinum (using up to 1536 cores), respectively, in a weak scaling
experiment. Here, the problem size per node is set to the data in
Table~\ref{tab:sizes}, and the throughput per node is compared. Ideal weak
scaling would correspond to coinciding curves for 1, 2, and 32 nodes. The performance
degradation is expected because MPI pack/unpack operations must be included,
bringing the observed performance closer to
the MPI-only performance. In this setting, the proposed Hermite-like basis is
able to maintain a much higher throughput than the nodal basis. This result confirms that large-scaling applications are
behaving more like the MPI-only case, where the proposed basis has a significant
advantage.

\section{Application to multigrid}\label{sec:multigrid}

For explicit time integration, the results in section~\ref{sec:experiments}
directly translate to application performance. For implicit solvers, the
matrix-vector product is embedded into
some iterative solver. In this section, we analyze the effect of the basis on
the throughput of a geometric multigrid solver using point-based smoothers and
cell-based block Jacobi smoothers.

\subsection{Algorithm}

Multigrid methods are highly efficient and scalable solvers for the linear systems
originating from elliptic partial differential equations. They combine simple
iterative schemes, effective in removing the high-frequency content, on
a hierarchy of coarser meshes. The hierarchy can be based on coarser meshes
(geometric multigrid), on lower polynomial degrees ($p$-multigrid), or on
algebraic coarsening based on the connectivity in the matrix (algebraic
multigrid) \cite{trottenberg01}. The most expensive
component of multigrid is usually the pre- and post-smoothing on the finest level.
In a massively parallel context,
the coarser level can contribute by the latency of the matrix-vector
products~\cite{Gholami16,Kronbichler18}, a cost that is ignored here because it is mostly
basis-agnostic.

For smoothing, a popular method in the context of matrix-free methods is the
Chebyshev iteration \cite{Adams03} around a simple scheme, like the inverse
of the matrix diagonal (point Jacobi) or some block-Jacobi/additive-Schwarz
method with block size equal to the number of unknowns per element. The
Chebyshev method is based on a three-term recurrence with iteration index $j$,
computing the solution $\boldsymbol u_{j+1}$ as follows,
\begin{equation}\label{eq:chebyshev}
\boldsymbol u_{j+1} = \boldsymbol u_{j} + \sigma_j (\boldsymbol u_{j}-\boldsymbol u_{j-1}) + \theta_j
P^{-1}\left(\boldsymbol b-A\boldsymbol u_{j}\right),
\end{equation}
where $\sigma_j, \theta_j$ are two scalar coefficients determined from Chebyshev
polynomials \cite{Adams03}, $P^{-1}$ is the inner preconditioner (e.g., point
Jacobi), $\boldsymbol b$ is the right hand side of the linear system, and $A$ the system
matrix. The iteration is started with an
initial guess $\boldsymbol u_{0}$ for the solution, e.g., from previous multigrid
iterations or the coarse grid correction, and $\sigma_0 = 0$. The main
computational expense in this method is the matrix-vector product, the
application of the preconditioner, and some vector updates.

\subsection{Efficiency with point-Jacobi/Chebyshev smoothing}

The simplest choice for the Chebyshev iteration is the point-Jacobi scheme with
$P=\text{diag}(A)$ in~\eqref{eq:chebyshev}. As the quality of the point-Jacobi
method strongly depends on the basis, this setup illustrates the properties
of the Hermite-like basis. We solve the 3D
Poisson equation on a cube $(-1,a_1)\times (-1,a_2)\times (-1,a_3)$ with
refinement selection and number of unknowns as in Table~\ref{tab:sizes}.
As before, we run the implementation for a general geometry and use a fixed
geometric coefficient for all quadrature points to highlight the cost of vector data access. The
lengths $a_i$ are chosen to ensure cuboidal elements depending on the
values of the mesh refinement $l_1, l_2, l_3$. If $l_1=l_2=l_3$, we set
$a_1=a_2=a_3=1$. If $l_1=l_2+1=l_3+1$, we set $a_1=3, a_2=a_3=1$.
Finally, if $l_1=l_2=l_3+1$, we set $a_1=a_2=3, a_3=1$.
The right hand side is chosen such that the
solution satisfies $u(x,y,z) = \sin(3\pi x)\sin(3\pi y)\sin(3\pi z)$.
Homogeneous Dirichlet conditions are set on all boundaries. We run a conjugate
gradient solver preconditioned by a geometric multigrid V-cycle
until the unpreconditioned $l_2$ norm of the residual,
$\|\boldsymbol b-A \boldsymbol u\|_2$, measured by residual estimate of the
conjugate gradient method,
has been reduced by $10^9$ compared to the initial residual $\|\boldsymbol b\|_2$ with
$\boldsymbol u^{(0)} = \boldsymbol 0$. The
multigrid V-cycle is completely done in single precision, whereas the
conjugate gradient solver runs in double precision. This gives essentially
identical accuracy as running the preconditioner in double precision but with
almost twice the throughput, see also \cite{Gropp00,Kronbichler19b}.
On each level, pre- and postsmoothing is done with a Chebyshev smoother of
degree 5. The relatively high degree is selected because it minimizes the time to solution with a
nodal Gauss--Lobatto basis over a range of polynomial degrees between 1 and 12
for the given mixed-precision setup.
The
Chebyshev parameters are selected to smoothen components in the range
$[0.06\tilde{\lambda}_{\max}, 1.2\tilde{\lambda}_{\max}]$ which is robust also
for mildly variable coefficients \cite{Kronbichler19b}. Here,
$\tilde{\lambda}_{\max}$ is an estimate of the largest eigenvalues found by 15
Lanczos iterations, starting with the vector $[-5.5, -4.5, \ldots, 4.5, 5.5,
-5.5, -4.5, \ldots]$. As a coarse-grid solver, we use the Chebyshev iteration
with the degree such that the a-priori error estimate of the Chebyshev
iteration guarantees a residual reduction by at least $10^{5}$ \cite{Varga09}.

\begin{table}
  \caption{3D Laplacian on a cube with MPI-only parallelization on 40 Xeon Gold
cores: Number of conjugate gradient iterations $n_9$ to reduce the $l_2$
residual by $10^9$ using geometric multigrid preconditioning with smoothing by a
Chebyshev iteration of point-Jacobi (first three column groups) as well as a
Chebyshev iteration with block-Jacobi via the fast diagonalization method (FDM)
\cite{lynch64} according to section~\ref{sec:iterations}. The degree of the
Chebyshev polynomial is 5 for both pre- and post-smoothing (i.e., 5
matrix-vector products).}
    \label{tab:mg_point_jacobi}
  {
	\footnotesize
    \setlength{\tabcolsep}{4pt}
	\strut\hfill
    \begin{tabular}{lcccccccccccc}
      \hline
      basis & \multicolumn{3}{l}{Hermite-like} & \multicolumn{3}{l}{nodal Gauss--Lob.} & \multicolumn{3}{l}{nodal Gauss} & \multicolumn{3}{l}{Hermite-like}\\
      $P^{-1}$ & \multicolumn{3}{l}{point Jacobi} & \multicolumn{3}{l}{point Jacobi} & \multicolumn{3}{l}{point Jacobi} & \multicolumn{3}{l}{block Jacobi: FDM}\\
      $p$ & $n_9$ & $\rho$ & MDoF/s & $n_9$ & $\rho$ & MDoF/s & $n_9$ & $\rho$ & MDoF/s & $n_9$ & $\rho$ & MDoF/s \\
      \hline
      1 & 8 & 0.058 & 26.3 & 8 & 0.058 & 26.4 & 8 & 0.060 & 29.2 & 7 & 0.042 & 29.7 \\
      2 & 26 & 0.45 & 8.99 & 7 & 0.036 & 30.0 & 9 & 0.085 & 24.7 & 7 & 0.034 & 33.1 \\
      3 & 29 & 0.48 & 8.77 & 7 & 0.047 & 33.2 & 11 & 0.13 & 23.0 & 7 & 0.044 & 36.1 \\
      4 & 16 & 0.27 & 16.6 & 7 & 0.038 & 32.2 & 12 & 0.17 & 20.2 & 7 & 0.035 & 37.6 \\
      5 & 13 & 0.20 & 21.1 & 8 & 0.057 & 28.1 & 15 & 0.24 & 15.5 & 7 & 0.051 & 38.1 \\
      6 & 12 & 0.17 & 22.4 & 7 & 0.051 & 30.0 & 16 & 0.27 & 13.8 & 7 & 0.046 & 37.3 \\
      7 & 13 & 0.20 & 19.3 & 9 & 0.085 & 19.3 & 20 & 0.34 & 9.13 & 8 & 0.068 & 30.3 \\
      8 & 14 & 0.22 & 18.9 & 9 & 0.085 & 18.9 & 21 & 0.35 & 8.42 & 8 & 0.076 & 31.8 \\
      9 & 15 & 0.23 & 17.0 & 10 & 0.12 & 15.3 & 24 & 0.42 & 6.58 & 9 & 0.086 & 27.7 \\
      10 & 15 & 0.24 & 16.9 & 10 & 0.11 & 15.2 & 25 & 0.43 & 6.24 & 9 & 0.091 & 27.2 \\
      11 & 16 & 0.26 & 14.7 & 12 & 0.16 & 12.0 & 29 & 0.48 & 5.12 & 10 & 0.12 & 23.0 \\
      12 & 16 & 0.27 & 11.6 & 12 & 0.16 & 9.83 & 29 & 0.48 & 3.73 & 10 & 0.11 & 18.9 \\
      \hline
    \end{tabular}
    \hfill\strut
  }
\end{table}

The results of the multigrid experiment for polynomial degrees $1\leq p \leq 12$
are presented in Table~\ref{tab:mg_point_jacobi}. Both
the number of iterations $n_9$ and the multigrid convergence rate
$\rho = \left(\nicefrac{\| \boldsymbol r_{n_9}\|_{2}}{\|\boldsymbol r_{0}\|_{2}}\right)^{1/n_9},$
involving the initial unpreconditioned residual $\boldsymbol r_0 = b$ and the residual $\boldsymbol r_{n_9}$
after $n_9$ iterations, are given. Furthermore, the table presents the
throughput of the multigrid solver in terms of million degrees of freedom solved
per second (MDoF/s), computed as the number of unknowns from
Table~\ref{tab:sizes} divided by the run time of the solver with an MPI-only
parallelization. The table compares the Hermite-like basis with a nodal
Gauss--Lobatto basis and a nodal Gauss basis for a point-Jacobi representation
of $P^{-1}$ in equation~\eqref{eq:chebyshev} as well as block-Jacobi variant
described in Section~\ref{sec:iterations} below. Among the three bases, the
nodal Gauss--Lobatto basis gives the best multigrid performance with point Jacobi and, as a
result, the highest solver throughput for $p\leq 6$. In other words, even though the
Hermite-like basis provides the fastest matrix-vector product according to the
right panel of Figure~\ref{fig:matvec_csl}, a simple
point-Jacobi smoother counteracts these benefits due to more multigrid iterations.
The nodal Gauss basis leads to relatively high iteration counts,
showing that optimality with respect to the mass matrix (which is diagonal in
the collocated Gauss case) does not translate to good performance for smoothing
with the
Laplacian. This also means that the Hermite basis could be somewhat improved
by using the additional degrees of freedom for orthogonality with respect to
the discrete Laplacian.

\subsection{Evaluation of $P^{-1}$}\label{sec:fdm}

A more general approach to remedy the disadvantage of the Hermite-like basis with point-Jacobi smoothing is a change of basis, which can recover
the behavior of any other DG basis. To achieve this,
we consider a preconditioner $P$ in the Chebyshev iteration \eqref{eq:chebyshev}
constructed as
\begin{equation}\label{eq:jacobi}
\left(P^{(K)}\right)^{-1} = \Big(T^\mathsf T \otimes T^\mathsf T \otimes
T^\mathsf T\Big)
\left(\hat{P}_{\text{d}}^\text{(K)}\right)^{-1} \Big(T \otimes T\otimes T\Big),
\end{equation}
where $\hat{P}_{\text{d}}^{(K)}$ is an element-wise diagonal matrix. The matrix
$T$ is some reference-cell 1D transformation matrix.
For example, a change into the nodal Gauss--Lobatto basis for
the sake of applying a diagonal operator $\hat{P}_{\text{d}}^{-1}$ allows the
Hermite-like basis to run with the same iteration counts as the nodal Gauss--Lobatto
basis (not shown) but with a considerably faster matrix-vector product.
For $p=7$, the nodal Gauss--Lobatto basis with a
point-Jacobi smoother solves 19.3
MDoF/s on 40 cores, whereas the Hermite-like basis including the transformation
\eqref{eq:jacobi} to the Gauss--Lobatto basis solves 26.7 MDoF/s.

\pgfplotstableread{
size timeTrafo timeDiag
16384     2.89733e-05 2.62777e-05
32768     2.41428e-05 2.27565e-05
65536     3.14168e-05 2.88138e-05
131072    3.48273e-05 3.14569e-05
262144    2.93464e-05 2.25855e-05
524288    3.64623e-05 2.51916e-05
1048576   6.08601e-05 2.89138e-05
2097152   7.60776e-05 5.39989e-05
4194304   0.000154743 0.000179451
8388608   0.000501386 0.000450909
16777216  0.00136468  0.00145859
33554432  0.00291979  0.00306269
67108864  0.00644327  0.00609875
134217728 0.0122768   0.0123495
268435456 0.0242952   0.0249275
}\tablePinvThree
\pgfplotstableread{
size      timeTrafo   timeDiag
21296     5.66883e-05 2.86941e-05
42592     5.25313e-05 2.11192e-05
85184     5.27336e-05 2.057e-05
170368    5.23477e-05 2.07019e-05
340736    5.96141e-05 2.82404e-05
681472    5.46489e-05 2.37414e-05
1362944   8.74909e-05 2.61458e-05
2725888   0.000164553 3.46808e-05
5451776   0.000305059 0.000241804
10903552  0.000859993 0.000898027
21807104  0.0019184   0.00176333
43614208  0.00388998  0.00382209
87228416  0.00799014  0.00783301
174456832 0.0159289   0.0158431
348913664 0.0320518   0.0321448
}\tablePinvTen
\pgfplotstableread{
degree size timeTrafo timeDiag
1  33554432 0.00345358 0.00317322
2  56623104 0.0056833  0.00524043
3  33554432 0.00306668 0.00304354
4  32768000 0.00282974 0.00280784
5  56623104 0.00541343 0.00533604
6  44957696 0.00394885 0.0039248
7  33554432 0.00291304 0.00305699
8  47775744 0.00437657 0.00422537
9  32768000 0.00278681 0.00296343
10 43614208 0.00392033 0.00382334
11 56623104 0.00520011 0.00498854
12 35995648 0.00341177 0.00314869
13 44957696 0.00448333 0.00403709
14 55296000 0.00575296 0.00496876
15 33554432 0.00388699 0.00290659
16 40247296 0.00441398 0.00359819
}\tablePinvDegree

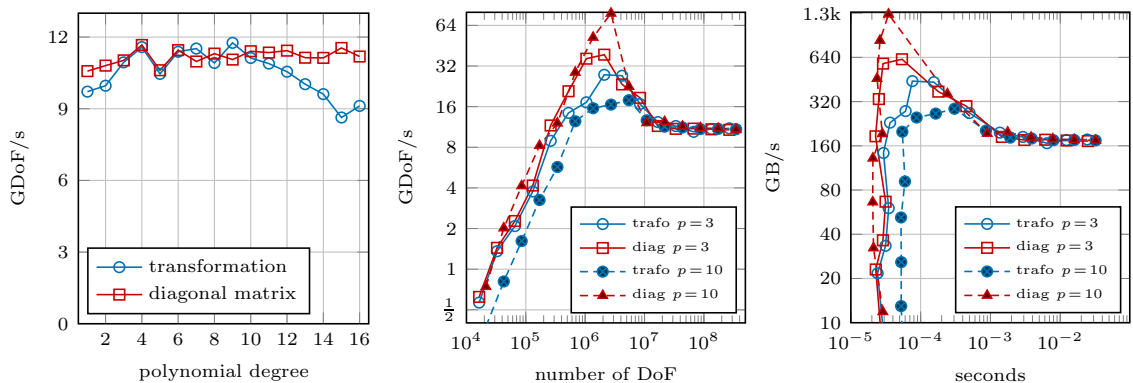
\begin{figure}
\strut\hfill
  \begin{tikzpicture}
    \begin{axis}[
      width=0.36\textwidth,
      height=0.38\textwidth,
      title style={font=\footnotesize},
      xlabel={polynomial degree},
      ylabel={GDoF/s},
      legend cell align={left},
      cycle list name=colorGPL,
      tick label style={font=\scriptsize},
      label style={font=\scriptsize},
      legend style={font=\scriptsize},
      legend pos=south west,
      grid,
      semithick,
      ymin=0,ymax=13,
      ytick={0,3,6,9,12},
      xmin=0.5,xmax=16.5,
      xtick={2,4,6,8,10,12,14,16},
      ]
          \addplot table [x expr={\thisrowno{0}}, y expr={1e-9*\thisrowno{1}/\thisrowno{2}/1}] {\tablePinvDegree};
          \addlegendentry{transformation};
          \addplot table [x expr={\thisrowno{0}}, y expr={1e-9*\thisrowno{1}/\thisrowno{3}/1}] {\tablePinvDegree};
          \addlegendentry{diagonal matrix};
    \end{axis}
  \end{tikzpicture}
  \hfill
  \begin{tikzpicture}
    \begin{loglogaxis}[
      width=0.35\textwidth,
      height=0.38\textwidth,
      title style={font=\footnotesize},
      xlabel={number of DoF},
      ylabel={GDoF/s},
      legend cell align={left},
      cycle list name=colorGPL,
      tick label style={font=\scriptsize},
      label style={font=\scriptsize},
      legend style={font=\tiny},
      legend pos=south east,
      grid,
      semithick,
      ymin=0.4,ymax=80,
      ytick={0.25,0.5,1,2,4,8,16,32,64},
      yticklabels={$\frac 14$,$\frac 12$,1,2,4,8,16,32,64},
      xmin=10000,xmax=5e8,
      xtick={1e4,1e5,1e6,1e7,1e8,1e9},
      ]
          \addplot table [x expr={\thisrowno{0}}, y expr={1e-9*\thisrowno{0}/\thisrowno{1}/1}] {\tablePinvThree};
          \addlegendentry{trafo $p\!=\!3$};
          \addplot table [x expr={\thisrowno{0}}, y expr={1e-9*\thisrowno{0}/\thisrowno{2}/1}] {\tablePinvThree};
          \addlegendentry{diag $p\!=\!3$};
          \addplot[gnuplot@darkblue,densely dashed,mark=otimes*,every mark/.append style={fill=gnuplot@darkblue!50!black,solid}] table [x expr={\thisrowno{0}}, y expr={1e-9*\thisrowno{0}/\thisrowno{1}/1}] {\tablePinvTen};
          \addlegendentry{trafo $p\!=\!10$};
          \addplot[gnuplot@red,densely dashed,mark=triangle*,every mark/.append style={fill=gnuplot@red!50!black,solid}] table [x expr={\thisrowno{0}}, y expr={1e-9*\thisrowno{0}/\thisrowno{2}/1}] {\tablePinvTen};
          \addlegendentry{diag $p\!=\!10$};
    \end{loglogaxis}
  \end{tikzpicture}
  \hfill
    \begin{tikzpicture}
    \begin{loglogaxis}[
      width=0.35\textwidth,
      height=0.38\textwidth,
      title style={font=\footnotesize},
      xlabel={seconds},
      ylabel={GB/s},
      legend cell align={left},
      cycle list name=colorGPL,
      tick label style={font=\scriptsize},
      label style={font=\scriptsize},
      legend style={font=\tiny},
      legend pos=south east,
      grid,
      semithick,
      ymin=10,ymax=1300,
      ytick={10,20,40,80,160,320,640,1280},
      yticklabels={10,20,40,80,160,320,640,1.3k},
      xmin=1e-5,xmax=0.099,
      xtick={1e-5,1e-4,1e-3,1e-2,1e-1},
      ]
          \addplot table [x expr={\thisrowno{1}}, y expr={16e-9*\thisrowno{0}/\thisrowno{1}/1}] {\tablePinvThree};
          \addlegendentry{trafo $p\!=\!3$};
          \addplot table [x expr={\thisrowno{2}}, y expr={16e-9*\thisrowno{0}/\thisrowno{2}/1}] {\tablePinvThree};
          \addlegendentry{diag $p\!=\!3$};
          \addplot[gnuplot@darkblue,densely dashed,mark=otimes*,every mark/.append style={fill=gnuplot@darkblue!50!black,solid}] table [x expr={\thisrowno{1}}, y expr={16e-9*\thisrowno{0}/\thisrowno{1}/1}] {\tablePinvTen};
          \addlegendentry{trafo $p\!=\!10$};
          \addplot[gnuplot@red,densely dashed,mark=triangle*,every mark/.append style={fill=gnuplot@red!50!black,solid}] table [x expr={\thisrowno{2}}, y expr={16e-9*\thisrowno{0}/\thisrowno{2}/1}] {\tablePinvTen};
          \addlegendentry{diag $p\!=\!10$};
    \end{loglogaxis}
  \end{tikzpicture}
  \hfill\strut
\caption{Analysis of evaluation of operator $P^{-1}$ with basis
  transformation~\eqref{eq:jacobi} versus a diagonal representation as
  function of the degree with 30--56 million DoF (left panel), as function of
  the size (middle panel), and as function of the evaluation time (right
  panel). The middle and right panels use the same data. All operations are
  done in single precision and with OpenMP on Xeon Gold.}
\label{fig:jacobi}
\end{figure}

To understand why the basis change leads to an increase in performance
despite additional arithmetic work, Figure~\ref{fig:jacobi} compares the
throughput of $P^{-1}$ with the basis change~\eqref{eq:jacobi} against the
application of a diagonal matrix (i.e., vector scaling) for experiments in
single precision. The throughput for large
sizes is essentially the same for both variants because they are both bound by
the memory
bandwidth of reading the vectors for $p\leq 11$, confirming previous results for the inverse mass matrix evaluation in
\cite{Fehn18b,Schoeder18} which is the same operation from an implementation point view \cite{Kronbichler16}.
In this experiment, $4\times 4=16$ Byte/DoF must be
accessed (read input \& diagonal, write and RFO for output), for 20--50 FLOP/DoF
for $1\leq p\leq 7$. For very high degrees, there is a slight degradation in throughput
for the transformation as the arithmetic work becomes notable. Experiments
that vary the size of the problem are shown in the middle and right panels of
Figure~\ref{fig:jacobi}. For an intermediate regime with enough
parallelism on the one hand and all data in the L2 and L3 caches on the
other hand, the diagonal preconditioner is faster
as expected. However, the transformed variant also runs much faster
when vectors are in cache, e.g., with 440 GB/s for $p=3$ compared to the saturated
performance of 180 GB/s. The in-cache case where
the diagonal matrix is advantageous is actually narrow: It only matters when the
time per iteration including the matrix-vector is around 0.5 milliseconds or
less. However, this regime is often beyond the strong scaling limit \cite{Fischer19}
where the network latency in the
matrix-vector products on
various levels is the dominant cost.

\subsection{Choice of inner preconditioner for Hermite-like basis}\label{sec:iterations}

The transformation from the Hermite-like basis into the nodal Gauss--Lobatto
basis for the purpose of smoothing is not the best one can do with
formula \eqref{eq:jacobi}. On a Cartesian (axis-aligned) mesh and with constant
coefficients, the tensor structure of shape functions and quadrature formula
propagates into the final cell matrix. Therefore, the discretization of the
scalar Laplacian admits an exact block-Jacobi preconditioner of the
form~\eqref{eq:jacobi}, where each block coincides with the cell-wise homogeneous
problem~\eqref{eq:discrete_laplace}.

For example, applying $\mathrm{det}(\mathcal J_{K}) = h^3$ and $\mathcal J_{K}^{-1} = \mathrm{diag}(h^{-1},h^{-1},h^{-1})$, where $h$ is the characteristic length of the Cartesian cell $K$, to the cell term~\eqref{eq:integral_laplace} with the ansatz function $u_h = \varphi_j$ it follows
\begin{equation}\label{eq:KDlaplace_cell}
(\nabla \varphi_i, \nabla \varphi_j)_{K}
=
\sum_{t=1}^3 \left(
\int_0^1 \frac{1}{h} \phi_{i_t}^\prime \phi_{j_t}^\prime \mathrm{d}\xi_t
\prod_{\tau=1, \tau \neq t}^3 \int_0^1 h \phi_{i_\tau} \phi_{j_\tau} \mathrm{d}\xi_\tau
\right).
\end{equation}
The associated stiffness matrix is a three-termed sum of Kronecker products of
the one-dimensional stiffness and mass matrices. Proceeding similarly with the
face integrals of the DG discretization, the matrix $A_K$ corresponding to
\eqref{eq:discrete_laplace} has the form
\begin{equation}\label{eq:KDlaplace_dg}
A^{(K)} = M \otimes M \otimes L_{\text{DG}} + M \otimes L_{\text{DG}} \otimes M + L_{\text{DG}} \otimes M \otimes M,
\end{equation}
where $L_{\text{DG}}$ represents the stiffness matrix of the 1D discretization
of~\eqref{eq:discrete_laplace} and $M$ the mass matrix on an interval with
length $h$. The generalized symmetric definite eigenproblem of the form
$L_{\textit{DG}} \boldsymbol z = \lambda M \boldsymbol z$ naturally leads to a basis transformation of
$\big(A^{(K)}\big)^{-1}$ of the form~\eqref{eq:jacobi} \cite{lottes05,lynch64}. Here,
the transformation matrix $T$ is the column-wise concatenation of generalized
eigenvectors $\boldsymbol z$ and the element-wise diagonal matrix is
\begin{equation}
\left(\hat{P}_{\text{d}}^{(K)}\right)^{-1} = \left(I \otimes I \otimes \Lambda + I \otimes \Lambda \otimes I + \Lambda \otimes I \otimes I\right)^{-1},
\end{equation}
where the diagonal matrix $\Lambda$ contains the generalized eigenvalues
$\lambda$ and $I$ is the identity matrix. Consequently, the tensor structure of
the DG disrectization of the scalar Laplacian admits an exact element-wise
preconditioner with diagonal form independent of the choice of the polynomial
basis. The inverse of $A^{(K)}$ in~\eqref{eq:KDlaplace_dg} is efficiently computed in terms of
the fast diagonalization method (FDM) \cite{lynch64}. We note that
\eqref{eq:KDlaplace_dg} needs to be slightly modified for elements at the
boundary due to the different weights resulting from the mirror principle. In
this work, we neglect this fact and use an approximate variant
\eqref{eq:KDlaplace_dg} with matrices from interior elements everywhere. For
further details on such block-based smoothers and possible extensions to
non-Cartesian meshes, see~\cite{WittePrep}.

The last three columns of Table~\ref{tab:mg_point_jacobi} report the iteration
counts and run time with an FDM-based representation in
equation~\eqref{eq:jacobi}. While the basis shows similar iteration counts as
the nodal Gauss--Lobatto case for $p\leq 4$, its multigrid convergence rates are
better for higher degrees. As explained above, combined with the basis change
in the smoother \eqref{eq:jacobi} the Hermite-like basis offers the best
performance in terms of time to solution.

\subsection{Efficiency of Chebyshev iteration}

Due to the high performance of the matrix-vector product, which is within a
factor of $1.7$ in throughput to simply copying
the involved vectors, the vector updates in the Chebyshev
update formula \eqref{eq:chebyshev} become critical. For optimal performance, they
must be merged with the other operations as much as possible. Also, the
isolated consideration of $P^{-1}$ as in Figure~\ref{fig:jacobi} is not enough.
In this section, we study two variants. The first variant, labeled
``separate'', computes $A\boldsymbol u^{(j)}$ and stores it in a temporary vector.
In a second loop through the data, the temporary vector is read again, $P^{-1}$ is applied,
and combined with the vectors. If $P^{-1}$ is represented by a vector (inverse
diagonal), this approach accesses 3 words of memory per unknown for the matrix-vector
product (including RFO) and 6 words during the preconditioner
evaluation and vector update phase, namely reading from $\boldsymbol b$, $A\boldsymbol u^{(j)}$,
$\hat{P}^{-1}_\text{d}$, $\boldsymbol u^{(j-1)}$, $\boldsymbol u^{(j)}$, as well as writing the result $\boldsymbol u^{(j+1)}$ back into the
storage location of $\boldsymbol u^{(j-1)}$ (this
avoids the RFO transfer). The second variant, labeled ``merged'', makes use of the fact that
the matrix-vector product computes the full result on an element within a single
sweep. Thus, after finishing the matrix-vector product on an element, we
immediately compute the residual, apply $P^{-1}$, and add the vector
contributions. In this variant, the temporary results as well as $\boldsymbol u^{(j)}$ are
still hot in caches, leading to a best-case data access of 5 words per unknown. All results in Table~\ref{tab:mg_point_jacobi} are based on the faster fully merged variant.

\pgfplotstableread{
degree size timeHermM  opsHerm timeGLM    opsGLM timeGM     opsGM timeHermS timeHermD  opsHermD
1  33554432 0.00913267 269     0.00931512 283    0.00810062 229   0.0125626 0.00904818 251
2  56623104 0.0154953  215.667 0.0175544  235    0.0159373  193   0.0192447 0.0134408  197.667
3  33554432 0.00763805 255     0.0078342  277    0.00712445 217   0.0108069 0.00697589 225
4  32768000 0.0065402  245     0.00719705 268.6  0.00668496 210.2 0.0094464 0.00626432 212.6
5  56623104 0.0111005  274.333 0.0124082  299    0.0114691  229   0.0158697 0.0107118  232.333
6  44957696 0.00899019 276.184 0.0104074  301.61 0.00965432 230.2 0.0127443 0.0085203  230.755
7  33554432 0.00716006 302     0.00840794 328    0.00761795 247   0.010014  0.00684805 248
8  47775744 0.00959643 309     0.0116474  335.44 0.0106756  251.6 0.0134174 0.00892072 251
9  32768000 0.00663308 333     0.0080811  359.8  0.00733941 267.4 0.0091508 0.00606959 267
10 43614208 0.00896905 342.785 0.0110127  369.87 0.0100699  273.9 0.0124919 0.00829538 272.421
11 56623104 0.0127239  365.667 0.0153509  393    0.0138815  289   0.017147  0.0113322  287.667
12 35995648 0.00844502 377.166 0.0103935  404.70 0.00943992 296.6 0.0115293 0.00768187 294.55
13 44957696 0.0108775  399.286 0.0133341  427    0.0121614  311.2 0.0152135 0.00992088 309.286
14 55296000 0.0137468  411.933 0.0172331  439.8  0.0158125  319.7 0.0193619 0.0125734  317.133
15 33554432 0.0105252  433.5   0.0132858  461.5  0.0119536  334   0.0151172 0.00990402 331.5
16 40247296 0.0120149  446.965 0.015689   475.08 0.0149186  342.9 0.0168074 0.0114691  340.024
}\tableChebyshev

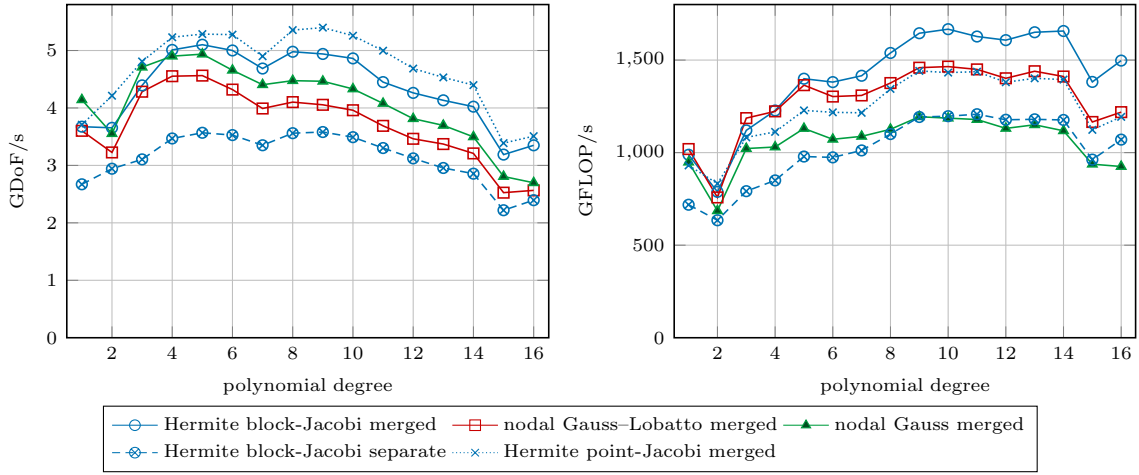
\begin{figure}
  \begin{tikzpicture}
    \begin{axis}[
      width=0.528\textwidth,
      height=0.4\textwidth,
      title style={font=\footnotesize},
      xlabel={polynomial degree},
      ylabel={GDoF/s},
      legend columns = 3,
      legend to name=legend:cslchebyshev,
      legend cell align={left},
      cycle list name=colorGPL,
      tick label style={font=\scriptsize},
      label style={font=\scriptsize},
      legend style={font=\scriptsize},
      grid,
      semithick,
      ymin=0,ymax=5.8,
      ytick={0,1,2,3,4,5,6},
      xmin=0.5,xmax=16.5,
      xtick={2,4,6,8,10,12,14,16},
      ]
      \addplot table [x expr={\thisrowno{0}}, y expr={1e-9*\thisrowno{1}/\thisrowno{2}/1}] {\tableChebyshev};
      \addlegendentry{Hermite block-Jacobi merged};
      \addplot table [x expr={\thisrowno{0}}, y expr={1e-9*\thisrowno{1}/\thisrowno{4}/1}] {\tableChebyshev};
      \addlegendentry{nodal Gauss--Lobatto merged};
      \addplot table [x expr={\thisrowno{0}}, y expr={1e-9*\thisrowno{1}/\thisrowno{6}/1}] {\tableChebyshev};
      \addlegendentry{nodal Gauss merged};
      \addplot[gnuplot@darkblue,densely dashed,mark=otimes,every mark/.append style={fill=gnuplot@darkblue!50!black,solid}] table [x expr={\thisrowno{0}}, y expr={1e-9*\thisrowno{1}/\thisrowno{8}/1}] {\tableChebyshev};
      \addlegendentry{Hermite block-Jacobi separate};
      \addplot[gnuplot@darkblue,densely dotted,mark=x,every mark/.append style={solid}] table [x expr={\thisrowno{0}}, y expr={1e-9*\thisrowno{1}/\thisrowno{9}/1}] {\tableChebyshev};
      \addlegendentry{Hermite point-Jacobi merged};
    \end{axis}
  \end{tikzpicture}
  \hfill
  \begin{tikzpicture}
    \begin{axis}[
      width=0.51\textwidth,
      height=0.4\textwidth,
      title style={font=\footnotesize},
      xlabel={polynomial degree},
      ylabel={GFLOP/s},
      legend pos={south west},
      legend cell align={left},
      cycle list name=colorGPL,
      tick label style={font=\scriptsize},
      label style={font=\scriptsize},
      legend style={font=\scriptsize},
      grid,
      semithick,
      ymin=0,ymax=1800,
      xmin=0.5,xmax=16.5,
      xtick={2,4,6,8,10,12,14,16},
      ]
      \addplot table [x expr={\thisrowno{0}}, y expr={\thisrowno{3}*1e-9*\thisrowno{1}/\thisrowno{2}/1}] {\tableChebyshev};
      \addplot table [x expr={\thisrowno{0}}, y expr={\thisrowno{5}*1e-9*\thisrowno{1}/\thisrowno{4}/1}] {\tableChebyshev};
      \addplot table [x expr={\thisrowno{0}}, y expr={\thisrowno{7}*1e-9*\thisrowno{1}/\thisrowno{6}/1}] {\tableChebyshev};
      \addplot[gnuplot@darkblue,densely dashed,mark=otimes,every mark/.append style={fill=gnuplot@darkblue!50!black,solid}] table [x expr={\thisrowno{0}}, y expr={\thisrowno{3}*1e-9*\thisrowno{1}/\thisrowno{8}/1}] {\tableChebyshev};
      \addplot[gnuplot@darkblue,densely dotted,mark=x,every mark/.append style={solid}] table [x expr={\thisrowno{0}}, y expr={\thisrowno{10}*1e-9*\thisrowno{1}/\thisrowno{9}/1}] {\tableChebyshev};
    \end{axis}
  \end{tikzpicture}
  \\
  \strut\hfill\pgfplotslegendfromname{legend:cslchebyshev}\hfill\strut
\caption{Throughput of one step within single-precision Chebyshev iteration
  ($j>0$) with the 3D Laplacian on an affine geometry for
  various bases with OpenMP parallelization on $2\times 20$ Xeon Gold
  cores. Both a block-Jacobi and diagonal preconditioner $P^{-1}$ are
  considered as well as separate matrix-vector product and fully merged iteration.}
\label{fig:chebyshev_csl}
\end{figure}

Figure~\ref{fig:chebyshev_csl} reports the throughput of one step in the
Chebyshev iteration \eqref{eq:chebyshev} when run for a generic index $j>0$
(i.e., $\boldsymbol u^{(j-1)}$ must be accessed) in single
precision. The fully merged variant is considerably faster due to the lower
vector access, reaching 5.0 GDoF/s per iteration for $p=8$ versus 3.6 GDoF/s when the
matrix-vector product runs separately. Note that the measured memory transfer of 120 GB/s
is around 30 GB/s less for the merged case and not fully utilizing the memory
interface.

When it comes to a point-Jacobi versus a block-Jacobi
variant according to \eqref{eq:jacobi}, Figure~\ref{fig:chebyshev_csl}
confirms that the basis change is almost for free. As opposed to
Figure~\ref{fig:jacobi} or when run with a slower separate matrix-vector product,
there is a small gap in throughput with the fastest implementation also for $p<12$
because the overall merged operator is not fully memory limited. Nonetheless,
the experiment is encouraging from a mathematical point of view because it
allows to select a basis where the diagonal is a good smoother: The
arithmetic operations for the basis change can mostly be hidden behind the
memory transfer in the Chebyshev loop also on an architecture with a
relatively low machine balance. The difference between the diagonal and
block-diagonal setup is even smaller on the Xeon Platinum with more
FLOP/Byte.

\section{Conclusions}
\label{sec:conclusions}

We have presented a Hermite-like basis that enables faster matrix-free
evaluation of symmetric interior penalty DG discretizations of the Laplacian
but at the same time yields the same results as a nodal Lagrangian basis
due to consistent integration. As opposed to
the higher-order continuity typically associated with Hermite polynomials, our
approach targets the fully discontinuous $L_2$-conforming case and is
motivated by the favorable data access. We have shown that the basis
significantly reduces the amount of vector data access on neighbors for higher
degree polynomials and especially with an MPI-only parallelization. Furthermore,
caches become more effective in holding neighbor data. We have shown
an 8--20\% increase in performance over nodal bases for an OpenMP
parallelization and up to 2$\times$ higher performance with an MPI
parallelization for polynomial degrees between 5 and 10. The basis is
specifically designed for modern hardware with high FLOP/Byte ratios, i.e., where
data access is expensive as compared to computations.

The proposed basis relies on a combination of the Hermite polynomials with nodal
polynomials based on roots of Jacobi polynomials to ensure a well-conditioned
interpolation. While the basis is not interpolatory,
it is constructed such that optimizations
of nodal codes like the even-odd decomposition are applicable.
Using the proposed basis with a point-Jacobi
smoother in a multigrid solver leads to worse behavior than with nodal
Lagrange polynomials based on Gauss--Lobatto points. We have therefore proposed
to combine fast operator evaluation in the Hermite-like basis with a basis
change into a more favorable basis for preconditioning. This basis could be the
nodal Gauss--Lobatto one, but more beneficial ones as exemplified by
the basis spanned by generalized eigenvectors from the fast diagonalization
method are also possible. The
change of basis happens in an otherwise memory-bandwidth limited algorithm and
can therefore be almost completely hidden behind the cost of memory transfer.

\appendix

\section*{Acknowledgments}
The authors thank Anian Fuchs and Guido Kanschat for
discussions about implementation aspects for the Hermite basis functions.

\section{Explicit formula of Hermite-like polynomials}\label{sec:appendix}
\noindent Basis functions for $p=2$:
\[
\phi_0(\xi) = (1-\xi)^2,\quad \phi_1(\xi) = 2\xi(1-\xi),  \quad \phi_2(\xi) = \xi^2.
\]
Basis functions for $p=3$:
\begin{align*}
&\phi_0(\xi) = -\frac 72 \left(\xi -\frac 27\right) (\xi-1)^2,\quad \phi_1(\xi) = \frac{11}{2} \xi (\xi-1)^2,\\
&\phi_2(\xi) = -\frac{11}{2} \xi^2(\xi-1),\quad \phi_3(\xi) = \frac 72 \xi^2\left(\xi-\frac 57\right).
\end{align*}
Basis functions for $p=4$:
\begin{align*}
&\phi_0(\xi) = 12 \left(\xi -\frac 16\right)\left(\xi -\frac 12\right) (\xi-1)^2, \quad
\phi_1(\xi) = -20 \xi \left(\xi -\frac 12\right) (\xi-1)^2,\\
&\phi_2(\xi) = 16 \xi^2 (\xi-1)^2,\\
&\phi_3(\xi) = -20 \xi^2\left(\xi -\frac 12\right)(\xi-1), \quad
\phi_4(\xi) = 12 \xi^2\left(\xi -\frac 12\right)\left(\xi-\frac 56\right).
\end{align*}
Basis functions for $p=5$:
\begin{align*}
&\phi_0(\xi) = -\frac{198}{5} \left(\xi -\frac 19\right)\left(\xi -\frac 12+\sqrt{\frac{1}{44}}\right)\left(\xi -\frac 12-\sqrt{\frac{1}{44}}\right) (\xi-1)^2, \\
&\phi_1(\xi) = \frac{1694}{25} \xi \left(\xi -\frac 12+\sqrt{\frac{1}{44}}\right)\left(\xi -\frac 12-\sqrt{\frac{1}{44}}\right) (\xi-1)^2, \\
&\phi_2(\xi) = -\frac{242\sqrt{44}}{25} \xi^2 \left(\xi -\frac 12-\sqrt{\frac{1}{44}}\right) (\xi-1)^2,\\
&\phi_3(\xi) = \frac{242\sqrt{44}}{25} \xi^2 \left(\xi -\frac 12+\sqrt{\frac{1}{44}}\right) (\xi-1)^2,\\
&\phi_4(\xi) = -\frac{1694}{25} \xi^2 \left(\xi -\frac 12+\sqrt{\frac{1}{44}}\right)\left(\xi -\frac 12-\sqrt{\frac{1}{44}}\right) (\xi-1),\\
&\phi_5(\xi) = \frac{198}{5} \xi^2\left(\xi -\frac 12+\sqrt{\frac{1}{44}}\right)\left(\xi -\frac 12-\sqrt{\frac{1}{44}}\right)\left(\xi-\frac 89\right).
\end{align*}
An algorithm for the polynomials for arbitrary order according to Section~\ref{sec:polynomials} is given in the deal.II library, version 9.0, class \texttt{HermiteLikeInterpolation},\\ \url{https://github.com/dealii/dealii/blob/dealii-9.0/source/base/polynomial.cc},\\
lines 1373--1494, retrieved on July 12, 2019.

{
\setlength{\itemsep}{0pt}
\setlength{\parsep}{0pt}
\setlength{\parskip}{0pt}

}

\end{document}